\documentclass[11pt,twoside]{article}
\usepackage[paper=a4paper,left=25mm,right=25mm,top=30mm,bottom=30mm]{geometry}
\usepackage{graphicx}
\usepackage{mathptmx}      
\usepackage{latexsym}
\usepackage{amssymb}
\usepackage{color}
\usepackage{natbib}
\usepackage{multirow}
\usepackage{rotating}
\usepackage{url}
\usepackage[breaklinks]{hyperref}
\hypersetup{
    colorlinks=true,
    linkcolor=black,
    citecolor=black,
    filecolor=black,
    urlcolor=black,
}
\newtheorem{hypothesis}{Hypothesis}
\usepackage{fancyhdr}
\usepackage{caption}
\begin{document}
\title{\LARGE \bfseries Test anxiety in mathematics among early undergraduate students in a British university in Malaysia}

\author{N. Karjanto$^{{\scriptsize\textmd{a,b}}}$\thanks{Corresponding author. Email: natanael@skku.edu} \ \ and \ S.~T.~Yong$^{\scriptsize \textmd{c}}$\\
$^{\scriptsize\textmd{a}}$\textsl{\small Department of Mathematics, School of Science and Technology, Nazarbayev University, Astana, Kazakhstan}\\
$^{\scriptsize\textmd{b}}$\textsl{\small Department of Mathematics, University College, Sungkyunkwan University,  Republic of Korea}\\
$^{\scriptsize\textmd{c}}$\textsl{\small Faculty of Science, The University of Nottingham Malaysia Campus, Semenyih, Malaysia}} 
\date{\footnotesize \textsl{(Received 29 April 2011; final version received 16 October 2012)}}


\maketitle

\begin{abstract}
\noindent
The level of test anxiety in mathematics subjects among early undergraduate students at The University of Nottingham Malaysia Campus is studied in this paper. The sample consists of 206 students taking several mathematics modules who completed the questionnaires on test anxiety just before they entered the venue for midterm exams. The sample data include the differences in the context of academic levels, gender groups and nationality backgrounds. The level of test anxiety in mathematics is measured using seven Likert questionnaire statements adapted from the Test Anxiety Inventory describing one's emotional feeling before the exam start. In general, the result shows that the students who had a lower score expectation were more anxious than those who had a higher score expectation, but that they obtained a better score than the expected score. In the context of academic levels, gender groups and nationality backgrounds, there were no significant correlations between the level of test anxiety and the students' academic performance. The effect size of the correlation values ranged from extremely small to moderate.\\

\end{abstract}

\section{Introduction}

In many academic courses, tests and examinations are usually used as primary tools to measure students' comprehension of the course material. As such, the validity and rigour of the examination process has received significant research attention. However, less attention has been paid to the exploration of the problems of the students who take the examinations (\citealp{Anderson95}). Since helping students to overcome test anxiety is one strategy to become a successful teacher (\citealp{Glasgow03}), test anxiety, particularly in connection to academic performance (AP) in mathematics, appears to be an issue worthy of research attention.

Test anxiety is a symptom or an emotion associated with stressful situations~(\citealp{Sarason84,Sarason90}), a feeling of apprehension and discomfort accompanied by cognitive difficulties during a test~(\citealp{Robinson09}). It is also defined as the set of phenomenological, physiological and behaviourial responses that accompany concern about possible negative consequences or failure in an exam or a similar evaluative evaluation~(\citealp{Zeidner98}). It involves a combination of physiological over-arousal, worry and dread about test performance and often interferes with normal learning and lowers test performance~(\citealp{Mandler52}).

In this article, the relationship among students' expectation, the level of test anxiety and their AP in several mathematics modules is investigated. The interest is to examine this relationship in the context of the study programme, gender group and local versus international differences. Furthermore, the correlation between students' expectation and their actual score in connection to the level of test anxiety in mathematics is also discussed. The research was conducted at the University of Nottingham Malaysia Campus (UNMC) during the period of January--June 2008. This research is geographically unique since its focus is on early undergraduate students of the University of Nottingham's branch campus in Malaysia. When this branch campus was opened in 2000, it became the first and the only British university in Malaysia and one of the first anywhere in the world. Although there are three other overseas universities opening their branches in the country, none of them is British university since all are Australian universities.\footnote{A couple years after this article is published, there are a couple of other UK universities in Malaysia, including Heriot-Watt in Putrajaya as well as The University of Southampton and Newcastle University Medicine, both in Educity Iskandar, Johor Bahru.}

This article is organized as follows. In the following three sections, some literature overview is presented. The difficulty that students face in mathematics and how this may give rise to test anxiety are discussed in Subsection~\ref{math}. The relationship between test anxiety and AP as well as the determination that test anxiety exists is the focus of Subsection~\ref{exist}. Some possible solutions or techniques to minimize test anxiety reported in the literature are discussed in Subsection~\ref{solution}. Next, six hypotheses are proposed in Section~\ref{hypo}. The first three hypotheses connect the triangular correlations among the students' expectation score in mathematics tests, the level of test anxiety during the test period and their actual AP in mathematics. The final three hypotheses consider the correlations between the students' AP and the level of test anxiety in mathematics in the context of academic levels, gender groups and nationality backgrounds. Section~\ref{method} provides essential information on the composition of the participants, instrument design, questionnaires measurement, data collection procedure and statistical analysis. Section~\ref{result} discusses some findings of this study. A correlation analysis is utilised to confirm whether the results obtained support the proposed hypotheses. The analysis includes the entire group analysis, the comparisons of the academic levels, gender groups and local versus international differences. Finally, Section~\ref{conclusion} presents some limitations, conclusions, remarks and future implications of this study.

\subsection{Mathematics and test anxiety} \label{math}

There is an exhaustive list of literature on test anxiety and mathematics anxiety and their relationship with AP. In this subsection, a literature overview pertaining to the difficulty that many students face in learning mathematics which in turn may give rise to test anxiety is presented. An excellent discussion on test anxiety and mathematics anxiety, including the reasons for the problems being so widespread, has been presented by \citeauthor{Hilton80a}~(\citeyear{Hilton80a, Hilton80b}). A brief review of the history and assessment of mathematics anxiety, its relationship with personal and educational consequences and its important impact on the measures of performance has been presented by~\citeauthor{Ashcraft09}~(\citeyear{Ashcraft09}). The authors conclude that mathematics anxiety causes a decline in AP when mathematics is performed under timed, high-stakes conditions, both in laboratory tests and in educational settings. This means that mathematics achievement and proficiency scores for math-anxious individuals are underestimates of true ability.

The Mathematics Anxiety Rating Scale (MARS) has been used to predict anxiety experiences and AP of Psychology and Mathematics students (\citealp{Morris78}). It has been found that the MARS scores for Psychology students are higher than those for the Mathematics students. These scores are useful predictors of both worry and emotionality and they are inversely related to the performance for the Psychology students. A strong inverse relationship has also been found between worry and performance for both groups and between emotionality and performance for the Psychology students.

The relationship of mathematics anxiety to situationally assessed test anxiety, mathematics performance, physiological arousal, and mathematics avoidance behaviour in undergraduate students has been investigated by \citeauthor{Dew84} (\citeyear{Dew84}). Their results indicate that mathematics anxiety measures are more highly rated to each other than to test anxiety. The functional similarities and differences in the cognitive processes involved in mathematics anxiety and test anxiety has been explored by \citeauthor{Hunsley87} (\citeyear{Hunsley87}). The results indicate that both mathematics anxiety and test anxiety account for unique variance in pre-examination appraisals, negative internal dialogue, post-examination appraisals, pre-examination and post-examination anxiety levels and several types of performance attributions. Only test anxiety, however, accounted for variance in subjects' actual examination performances. The implications of a cognitive theory for mathematics anxiety have also been addressed in that paper.

The relationship among test anxiety, mathematics anxiety and teacher feedback among university students has been discussed by \citeauthor{Green90}~(\citeyear{Green90}). The author discovered that test anxiety has a greater effect on the mathematics achievement of remedial mathematics students than either mathematics anxiety or teacher comments. In addition, teacher feedback in the form of free comments and specified comments is more facilitative of test performance than no comments at all.

It is generally understood that that students' different levels of test anxiety are related to different subject areas. A variety of techniques should be used in helping students to overcome their test anxiety. One study compared 196 first-year college students' self-report of test anxiety in mathematics, physical sciences, English and social studies~(\citealp{Everson93}). The results indicate that subjects in the rank order of causing the most test anxiety to least test anxiety are physical sciences, mathematics, English and social studies.

Some notes on mathematics anxiety and test anxiety among diverse populations of students based on the author's experience, including cause, effect, pros and cons, can be found in \citeauthor{Marikyan09}~(\citeyear{Marikyan09}). Another study from Nigeria shows that many students have fear of mathematics and mathematics tests. The findings also reveal that the majority of the students do not know how to study for mathematics tests~(\citealp{Olatunde09}).

\subsection{Test anxiety and AP} \label{exist}

This subsection discusses a literature overview on the relationship between test anxiety and AP. The relationship of cognitive (worry) and emotional (emotionality) components of test anxiety with physiological arousal and AP among both high school and college students has been studied by \citeauthor{Morris70}~(\citeyear{Morris70}). The authors discovered two important aspects: worry is more highly negatively related to examination grades than is emotionality or pulse rate and worry is more highly negatively related to expectancy than is emotionality.

Study behaviors and their relationships with test anxiety and AP have been discussed by \citeauthor{Allen72}~(\citeyear{Allen72}). The relationship of test anxiety and AP when controlled for intelligence has been studied by \citeauthor{Boor72}~(\citeyear{Boor72}). The effects of the knowledge of ability test results on AP and test anxiety have been explained by \citeauthor{McMahon73}~(\citeyear{McMahon73}). The author hypothesized that the students with detailed knowledge would obtain higher grade point averages (GPAs) and have a lower level of test anxiety with no significant differences in GPAs among the groups.

The relationship of test anxiety with AP in college students, the differences in study-related behaviours between high and low test-anxious students and the differential effectiveness of study-related behaviours for both groups have been investigated by \citeauthor{Culler80}~(\citeyear{Culler80}). The results demonstrate a significant reduction in GPAs associated with test anxiety and high test-anxious students have also been found to have poorer study skills. For high test-anxious students, the quality of study habits and amount of study time were positively related to AP, whereas missing classes and delaying examinations were inversely related to AP.

A cross-cultural study of the relationship among AP, test anxiety, intelligence and sex has been conducted by \citeauthor{vanderPloeg84}~(\citeyear{vanderPloeg84}). The authors observed that boys and girls with lower intelligence achieved less and were less influenced by the impairing effects of test anxiety. The impact of test anxiety on test performance and the cognitive appraisals of test-anxious undergraduate students who enrolled in a Statistics course have been explored by \citeauthor{Hunsley85}~(\citeyear{Hunsley85}). The results indicate that test anxiety is related to poor test performance both early and late in the term. When state anxiety levels were controlled for, the test anxiety-test performance relation was apparent only during the later stages of the course. The pattern of students' anxiety and appraisals suggests that test-anxious students experience most doubt and concern early in the term.

Another study shows that test anxiety causes poor AP and it relates inversely to students' self-esteem and directly to their fears of negative evaluation, defensiveness and other forms of anxiety~(\citealp{Hembree88}). The relationship between test anxiety and AP in the context of the time of testing, that is, pre-test versus post-test, has been studied by \citeauthor{Zeidner91}~(\citeyear{Zeidner91}). The author verified the hypothesis that anxiety will be more strongly related to test performance at post-test than at pre-test. Whereas under pre-test anxiety measurement conditions, only a negligible relationship was observed between performance and test anxiety, under post-test measurement conditions, the latter two variables were observed to be moderately correlated.

Cognition, study habits, test anxiety and AP have been studied by \citeauthor{Kleijn94}~(\citeyear{Kleijn94}).
Another relationship between cognitive test anxiety and AP has been discussed by \citeauthor{Cassady02}~(\citeyear{Cassady02}). The authors found that a higher level of cognitive test anxiety was associated with significantly lower test scores. Their results support the conclusion that cognitive test anxiety exerts a significant, stable and negative impact on AP measures. The relationships between test anxiety experienced by students at the time of the final examination and students' performance expectations (PEs), actual performances and level of preparation for the final examination have been examined by \citeauthor{Burns04}~(\citeyear{Burns04}). The results support the relationship hypothesized between test anxiety and PEs at the time of the final examination and that the level of anxiety is greatly affected by the students' expectation. Furthermore, a study of the relationship between test anxiety and AP in undergraduate and graduate students has been carried out by \citeauthor{Chapel05}~(\citeyear{Chapel05}). The authors found a significant but small inverse relationship between test anxiety and GPAs in both groups. The current issues in test anxiety theory, assessment, research and intervention, including the relationship between test anxiety and AP, have been discussed by \citeauthor{Zeidner07}~(\citeyear{Zeidner07}).

A systematic inquiry upon the levels of test anxiety through behavioral and physiological procedures before and after a mathematics test has been done by \citeauthor{Martinez08}~(\citeyear{Martinez08}). It was found that academic general scores were inverse to the behavioural anxiety level, but the best mathematics scores were associated with the middle levels of behavioural anxiety. A study of the level of cognitive test anxiety of selected undergraduate students in Nigeria has been conducted by \citeauthor{Faleye10}~(\citeyear{Faleye10}). The author discovered that although the level of cognitive test anxiety negatively affects the performance level, the sex differences do not lead to the corresponding differences in both test anxiety and AP.

An investigation of possible effects on skills training of emotions confronting on a decrement of stress and test anxiety and on an improvement of AP of female high school students has been done by \citeauthor{Zohreh10}~(\citeyear{Zohreh10}). The relationships among achievement goal orientation, test anxiety, self-efficacy, metacognition and AP have been investigated in a recent Master's degree dissertation (\citealp{Lindsay10}). Another finding shows that students who normally experience a high level of test anxiety in the classroom have reduced test anxiety when taking online examinations, while the reverse is true for those who have a low level of anxiety in the classroom. Furthermore, the relationship between test anxiety and examination performance is weaker in an online setting than in the classroom~(\citealp{Stowell10}). The authors recommend that instructors evaluate the potential impact of these findings when considering offering examinations online. Implementation intentions and test anxiety in connection to shielding AP from distraction have been discussed by \citeauthor{Parks-Stamm10}~(\citeyear{Parks-Stamm10}). The results show that as test anxiety increases, the effectiveness of temptation-inhibiting implementation intentions increases too.

A relationship between a cognitive test anxiety component of `worry' and self-doubt is discussed by \citeauthor{Schwarzer96} (\citeyear{Schwarzer96}).
A review on findings from several studies on cross-cultural differences in academic achievement, anxiety and self-doubt, with a focus on comparisons between Confucian Asian and European regions, has been presented by \citeauthor{Stankov10}~(\citeyear{Stankov10}). The author reasons that the relatively unforgiving attitude among Confucian Asian students, coupled with the belief that effort rather than ability is the primary source of success, may be able to explain both high achievement and high anxiety of the group. In addition to that, the correlation among test anxiety, perfectionism, goal orientation and AP  has been discussed recently by \citeauthor{Eum11}~(\citeyear{Eum11}), who found that cognitive test anxiety is inversely associated with both performance indicators and positively associated with maladaptive perfectionism and avoidance goal orientations. Their results suggested that students who are highly test anxious are likely to be women who endorse avoidance goal orientations and are maladaptively perfectionists.

\subsection{Minimisation of test anxiety} \label{solution}

This subsection presents a literature overview on some possible technique to reduce the level of test anxiety. For instance, a student who experienced test anxiety was successfully treated by cue-controlled relaxation methods~(\citealp{Russel73}). The procedure involved training in deep-muscle relaxation and pairing of breath exhalations while relaxed with a self-produced cue word of calm. A literature review in behavioural approaches to the treatment of test anxiety has also been presented by \citeauthor{Spielberger76}~(\citeyear{Spielberger76}). The authors concluded that the successful reduction of anxiety is related to improving deficient study habits that interfere with adequate test preparation. Another research shows that contrary to prior perceptions, improved test performance and GPAs consistently accompany test anxiety reduction~(\citealp{Hembree88}). The same author also discussed a variety of treatments which are effective in reducing mathematics anxiety~(\citealp{Hembree90}), including classroom interventions and out-of-class psychological treatments.

It has been shown that guidance-oriented periodic testing has a positive effect on the performance of the students majoring in Engineering in Nigeria. However, none of the treatments has an effect on the mathematics-test anxiety of the students~(\citealp{Adedayo96}). An interesting technique to reduce anxiety and improve AP has been studied by \citeauthor{Beauchemin08}~(\citeyear{Beauchemin08}). The authors suggested the application of a mindful meditation and relaxation training. Mathematics anxiety can be treated either with direct or indirect interventions such as relaxation therapy, teaching style and cooperative learning~(\citealp{HellumA10}). Some tips on conquering mathematics anxiety can be found in \citeauthor{Arem10}~(\citeyear{Arem10}).

In the following section, a number of proposed hypotheses related to the study and some of the related references are considered.

\section{Hypotheses} \label{hypo}

The aim of this study is to investigate the relationships between students' reported levels of anxiety at the time of a test, within 20-30 minutes of starting their mathematics midterm examinations, and their AP. The following hypotheses are proposed.

At the time of test, the students have already been exposed to other type of assessments, for example, assignments and quizzes. Those who have performed well on the previous assessments are expected to be less apprehensive about the test. These students are believed to possess greater confidence and thus be less anxious about the test. Similarly, the students who performed poorly will relate such failures to personal failure. Those students have a lower expectation and exhibit a higher level of test anxiety~(\citealp{Burns04,Cassady02}). This leads to the first hypothesis.
\begin{hypothesis}
  There is a significantly negative correlation between the students' expectation of AP in the mathematics test and the level of test anxiety in mathematics during the test period. \label{hypothesis1}
\end{hypothesis}

Past research suggests that anxiety can impair AP among the students. In other words, it can be expected that individuals experiencing a higher level of test anxiety at the time of test will perform more poorly than the students who experience a lower degree of anxiety. Test anxiety in general is expected to have a negative effect on AP (\citealp{Smith64,McMahon73,Culler80,Hunsley85,Yee87,Hembree88,Horn89,Cassady02,Burns04,Chapel05}). This leads to the second hypothesis.
\begin{hypothesis}
  There is a significantly negative correlation between the students' level of test anxiety in mathematics during the test period and their AP in mathematics.
  \label{hypothesis2}
\end{hypothesis}

As the third side of the correlation triangle, it would be interesting to test another hypothesis that connects the students' expectations and their obtained actual scores. It is reasonable to assume that the students with a high expectation would also obtain a high score in their test. In particular, the relationship between students' expectations and their actual performances has been investigated by \citeauthor{Burns04}~(\citeyear{Burns04}). Another study shows that students' expectations and their AP did not correlate significantly (\citealp{Alexitch88}). This leads us to the third hypothesis.
\begin{hypothesis}
  There is a significantly positive correlation between the students' expectation on AP in the mathematics test and the their actual AP in mathematics.
  \label{hypothesis3}
\end{hypothesis}

Higher level students may have experienced or developed techniques to deal with test anxiety that freshmen may not have acquired. They may be more ready to self-identify and seek assistance~(\citealp{Ross06}). A study conducted in the USA showed that the levels of test anxiety were higher among high school students than among college students~(\citealp{Sud91}). This leads to the fourth hypothesis.
\begin{hypothesis}
  There is a significant difference in AP between the students enrolled in the pre-undergraduate (Foundation) programme and the ones enrolled in the Undergraduate programme when dealing with test anxiety in mathematics.
  \label{hypothesis4}
\end{hypothesis}

It has been found that female students experience higher levels of test anxiety than males irrespective of cultural background~(\citealp{vanderPloeg84,Sharma90,Kleijn94,Karimi09}). However, in this study, an examination goes further whether there is a significant difference in AP between the gender groups when dealing with test anxiety. This leads to the fifth hypothesis.
\begin{hypothesis}
  There is a significant difference in AP between male and female students when dealing with test anxiety in mathematics.
  \label{hypothesis5}
\end{hypothesis}

It is generally presumed that overseas students have a higher level of test anxiety than the local ones. It is also quite common that when people are exposed to different environments, particularly where habit and culture are different from the ones at home, they might get more anxious when under pressure. For instance, it has been discovered that compared with Australian students, other students mainly from Asian regions of Singapore, Malaysia and Hong Kong have significantly greater difficulties in adjusting to academic requirements~(\citealp{Burns91}). Furthermore, the author also reported significantly higher level of anxiety among these Asian students. There are also several other comparisons in terms of anxiety and AP across cultural differences, for instance, between native English-speaking children and bilingual Latino or Spanish dominant children~(\citealp{Willig83}), the anxiety of mathematics learning between Tanzanian and Malaysian students~(\citealp{Mohamed10}) and the level of anxiety and AP among Confucian Asian and European students~(\citealp{Stankov10}). This leads to the final hypothesis.
\begin{hypothesis}
  There is a significant difference in AP between local Malaysian and international students when dealing with test anxiety in mathematics.
  \label{hypothesis6}
\end{hypothesis}

The research methodology, which includes the composition of participants, instrument design and measurement, procedure and statistical analysis, is explained in the following section.

\section{Method} \label{method}

\subsection{Participant}

The participants of this research were the students enrolled at the Faculty of Engineering at the UNMC who took a number of mathematics modules as part of their curriculum. Out of an estimated number of around 250 distributed questionnaires, a sample of 206 valid questionnaires was returned. All these students can be classified according to the academic levels, gender groups and nationality backgrounds, that is, local Malaysians or international students. Although the students from overseas belong to a variety of nationality backgrounds, all of them are simply classified as international students. Based upon the sample data of academic levels, 156 students are enrolled in the Foundation programme (75.73\%) and 50 students are enrolled in the Undergraduate programme (24.27\%). Based upon the sample data on gender groups, 151 are male students (73.30\%) and 55 are female students (26.70\%). Based upon the sample data on nationality backgrounds, there are 148 local Malaysian students (71.84\%) and 58 international students from overseas (28.16\%). Table~\ref{data1} presents the composition of polled students according to their backgrounds.
\begin{table}[h]
{\renewcommand{\tabcolsep}{15pt}
\begin{center}
{\footnotesize
\caption{The composition of polled students.} \label{data1}
\begin{tabular}{lrcrrc}
  \hline \hline
         & Foundation & Undergraduate & Total & Local & International \\ \hline
  Male   & 102 & 49 & 151 & 105 & 46\\
  Female \qquad & 38  & 17 & 55  & 43  & 12 \\
  Total  & 140 & 66 & 206 & 148 & 58\\
  \hline \hline
\end{tabular}
}
\end{center}
}
\end{table}

\subsection{Instrument design and Measurement}

A two-page questionnaire was distributed to as many students as possible at several occasions before they entered the venues for the midterm examination in mathematics. It is estimated that the total number of students who have been given forms is around 250. Page~1 of the questionnaire consisted of a brief description of the questionnaire, including its main objective, a request to fill it in and the instruction as how to do so. Page~2 of the questionnaire consisted of two parts. Part~A collected the background information of the students. This included student identity number, programme enrolled, gender, nationality and the expected score (out of 100) for that particular test. The purpose of collecting the student identity number was to determine the actual test scores for each student after all the examination scripts had been marked. Regarding nationality, the options were either local (Malaysian) or international students. For the latter case, the students were also asked to specify their nationality.

Part~B of the questionnaire focused on the measurement of the students' level of test anxiety in mathematics. This part consists of seven Likert questionnaire statements, or simply, Likert items used to measure their level of test anxiety.\footnote{Likert items and Likert scales are types of psychometric items and scales frequently used in psychology questionnaires. A~Likert item is a statement which the respondent is asked to evaluate according to any kind of subjective or objective criteria where generally the level of agreement or disagreement is measured. A~Likert scale is the sum of responses on several Likert items. The scale was developed by and named after an organizational psychologist Rensis Likert~(\citeyear{Likert32}).} These Likert statements are based on the existing measurement that has been conducted in other studies, for instance in that of \citeauthor{Hong99}~(\citeyear{Hong99}). Similar to her research, the implemented questionnaire was derived from a modified version of the Test Anxiety Inventory (TAI)~(\citealp{Spielberger80}). This TAI is one of the most widely used instruments for measuring test anxiety for both high school and university students~(\citealp{Chapel05}). It is a self-report instrument that consists of 20 statements to measure students' emotional state and worry around the time of their tests, for which the respondents are asked to report how often they experience symptoms before, during and after taking the tests. Cronbach's $\alpha$ for this instrument is 0.92.\footnote{Cronbach's $\alpha$ (alpha) is a statistic that has an important use as a measure of the reliability of a psychometric instrument. It was first named as alpha by Lee Cronbach (\citeyear{Cronbach51}), as he had intended to continue with further instruments. Cronbach's $\alpha$ will generally increase when the correlations between the items increase. For this reason, the coefficient is also called the `internal consistency' or the `internal consistency reliability' of the test. As a rule of thumb, a requirement of $\alpha \geq 0.70$ is considered reliable.} The main reason for adopting the seven items from TAI is brevity so that the students could fill in the questionnaires within a short period of time. The questionnaire employed in this research, however, was designed to address the state of anxiety at one occasion only, that is, immediately before the examination.

When filling in Part~B, the students were asked to choose only one option that suited their condition at that moment the best. For each of the corresponding statement, the options were given by four-point Likert items ranging from 1 until 4, where each number denotes strongly disagree, disagree, agree and strongly agree, respectively. Since there is no neutral option, a `forced choice' method is thus implemented where the students are forced to choose either favouring the statement or otherwise. It has been shown that the overall difference between a four-point and a five-point Likert scale is negligible~(\citealp{Armstrong87}). A~reliability of this measurement has been checked and the result shows an acceptable value of Cronbach's $\alpha$ of $0.84$. The sample questionnaire is given in the appendix of this article.

\subsection{Procedure}

It is a custom at the UNMC that the students are not allowed to enter the examination venue prior to an instruction from the invigilators. Generally, many students would congregate in the corridor near the examination venues while they are checking their notes or doing some preparation for the examination. Between 10 and 20 minutes before the start of the examination, there was an opportunity to distribute the questionnaires. During this period, the students were approached and were asked to fill in the forms voluntarily. Around 5-10 minutes before the start of the examination, after the students were allowed to enter the venue of the examination, there was no further request to fill in the questionnaires.

Since many of them were already under pressure, it was emphasized that filling the questionnaire will only take around one to two minutes. Although many were willing to volunteer, some of them refused to do so and they were not forced to fill in the questionnaire. For those who were not willing to volunteer, it has been concluded that most likely the students did not prefer to be bothered and simply wanted to use the remaining time to concentrate on their final preparation for the test. Another possible reason for this might be that the questionnaire sheet does not contain a consent form.

\subsection{Statistical analysis}

Several simple statistical analyses were carried out in this research: descriptive Pearson correlation and $t$-test. Data analysis was performed using Microsoft Excel and SPSS 16.0 software. The histogram and scatter plot figures were generated using Matlab R2010b.

\section{Result} \label{result}

\subsection{General overview}

It is understood that generally students expect higher scores in the examinations than what they could achieve, which is also confirmed in this study. The following are the results of the expected score and the actual AP averages for different categories. A possible maximum score is 100 points. It is noticed that the average of the expected scores for all students is 76.44, while the average of the actual scores is 58.26. Based upon academic levels, the expected score averages for the students enrolled in the Foundation and Undergraduate programmes are 75.90 and 78.12, respectively. The actual score averages for the Foundation and Undergraduate students are 55.58 and 66.60, respectively. Based upon gender groups, the expected score averages for the male and females students are 76.75 and 75.58, respectively. The actual score averages for the male and female students are 55.21 and 66.63, respectively. Based upon nationality backgrounds, the expected score averages for the local and international students are 76.98 and 75.07, respectively. The actual score averages for the local and international students are 62.72 and 46.86, respectively. The average expected and actual scores for all students, students enrolled in the Foundation and Undergraduate programmes, male and female students and local and international students are listed in Table~\ref{average1}. 

\begin{table}[h]
\begin{center}
{\footnotesize
  \caption{The averages of the expected score and the actual outcome for all students, the Foundation and Undergraduate students, the male and female students as well as the local and international students.} \label{average1}
    \begin{tabular}{lccccccc}
       \hline \hline
                      & \quad All \quad \qquad & Foundation & Undergraduate &
                      \quad Male \quad \quad & Female \quad \quad & Local & International \\ \hline \hline
Expected score &  76.44   &  75.90   &  78.12 ($p = 0.493$) &  76.75    & 75.58 ($p = 0.709$) & 76.98 & 75.07\\
Actual score   &  58.26   &  55.58   &  66.60 ($p = 0.495$) &  55.21    & 66.63 ($p = 0.009$)        & 62.72 & 46.86\\ \hline
    \end{tabular}
}
\end{center}
\end{table}

\begin{figure}[h]
\begin{center}
  \includegraphics[width=0.25\textheight]{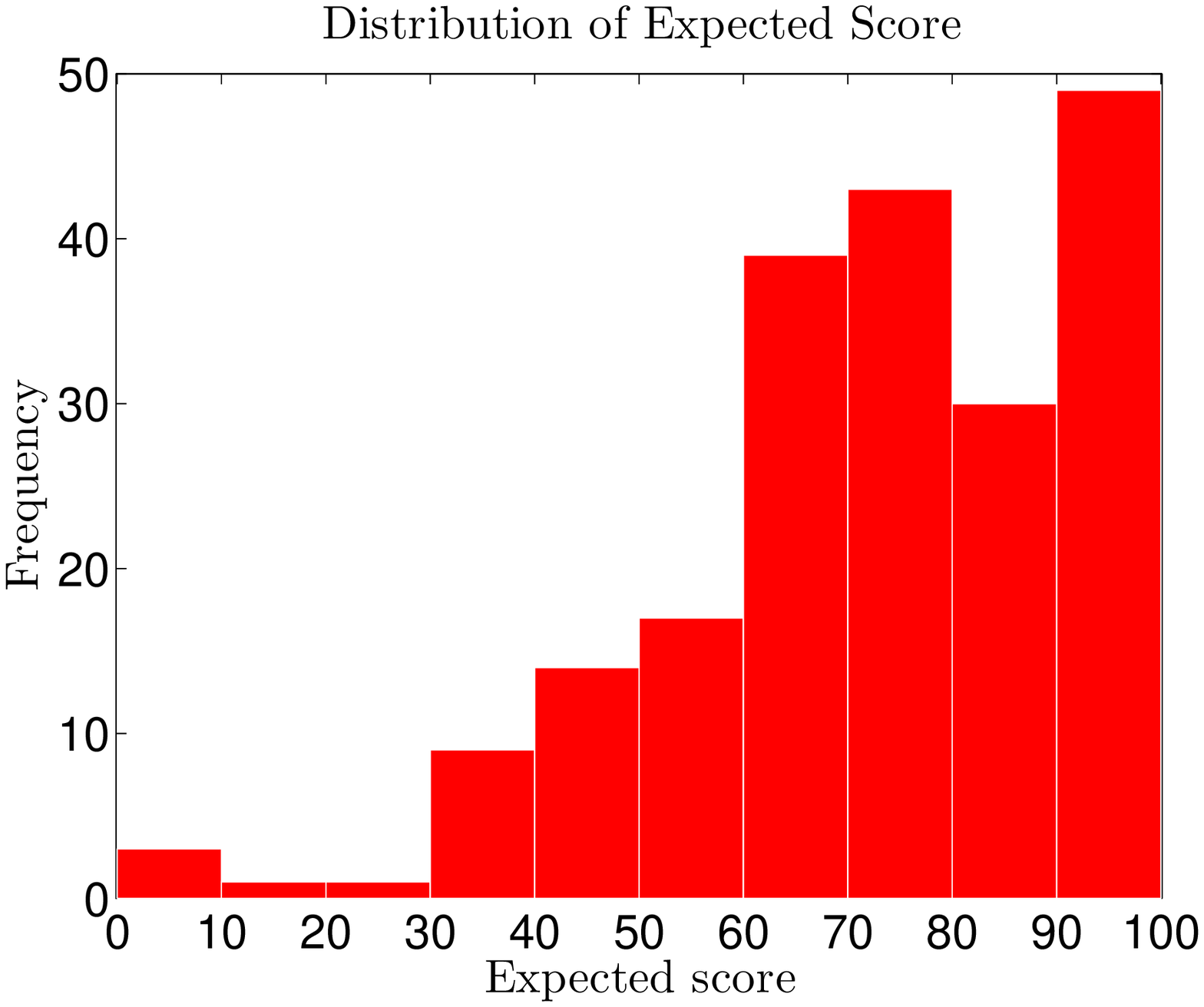} \hspace*{2cm}
  \includegraphics[width=0.25\textheight]{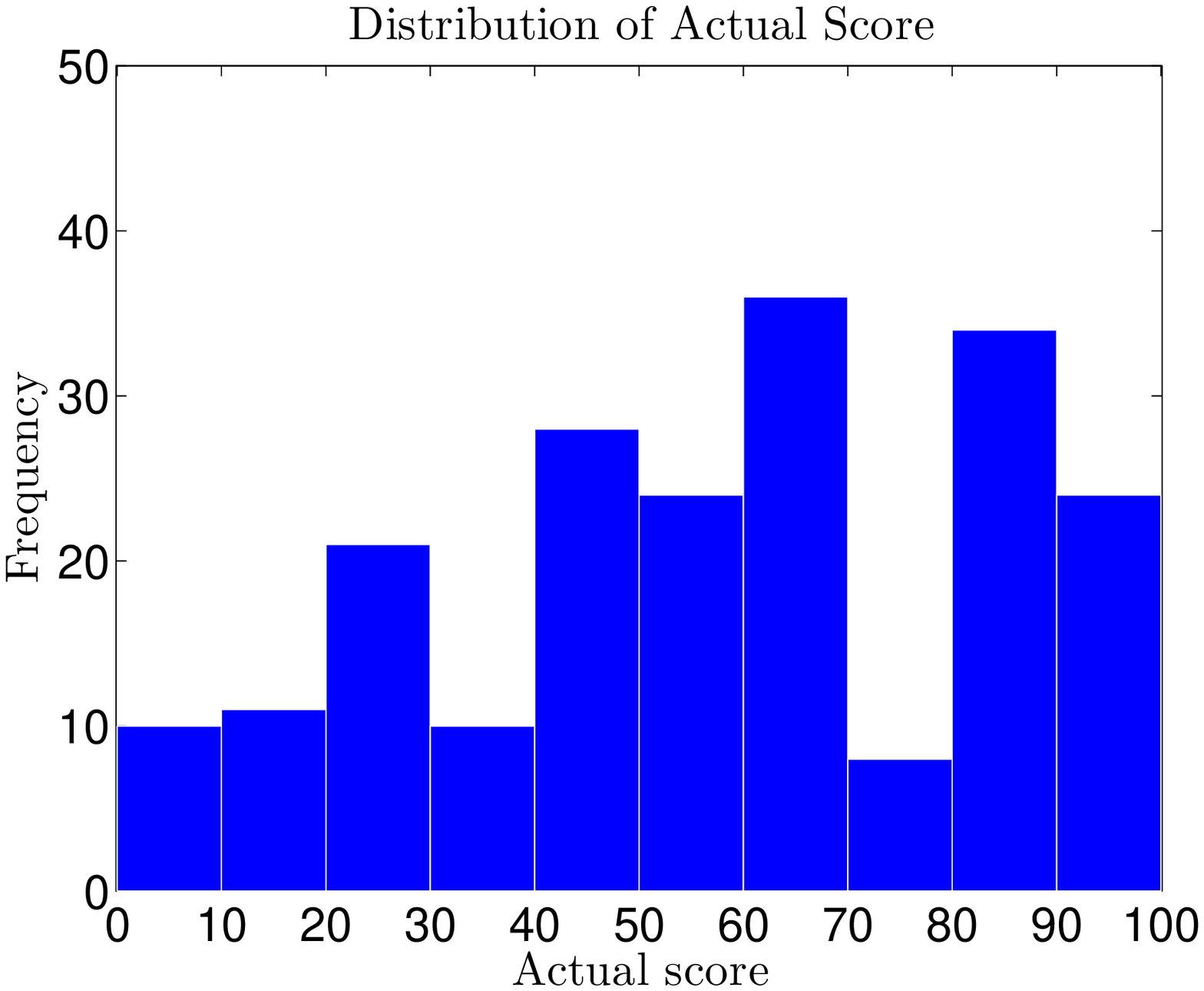} \\
  \caption{Histograms based on the collected sample data of the expected score (red histogram, left panel) and the actual AP (blue histogram, right panel). From these histograms, it can be observed that the expected score data show a non-normal left-skewed type of distribution where more students tend to expect high scores in the test, while the actual score data show a multimodal (plateau) type of distribution where more data are distributed to each score interval.}
  \label{Histo1}
\end{center}
\end{figure}

\begin{figure}[h]
\begin{center}
  \includegraphics[width=0.42\textheight]{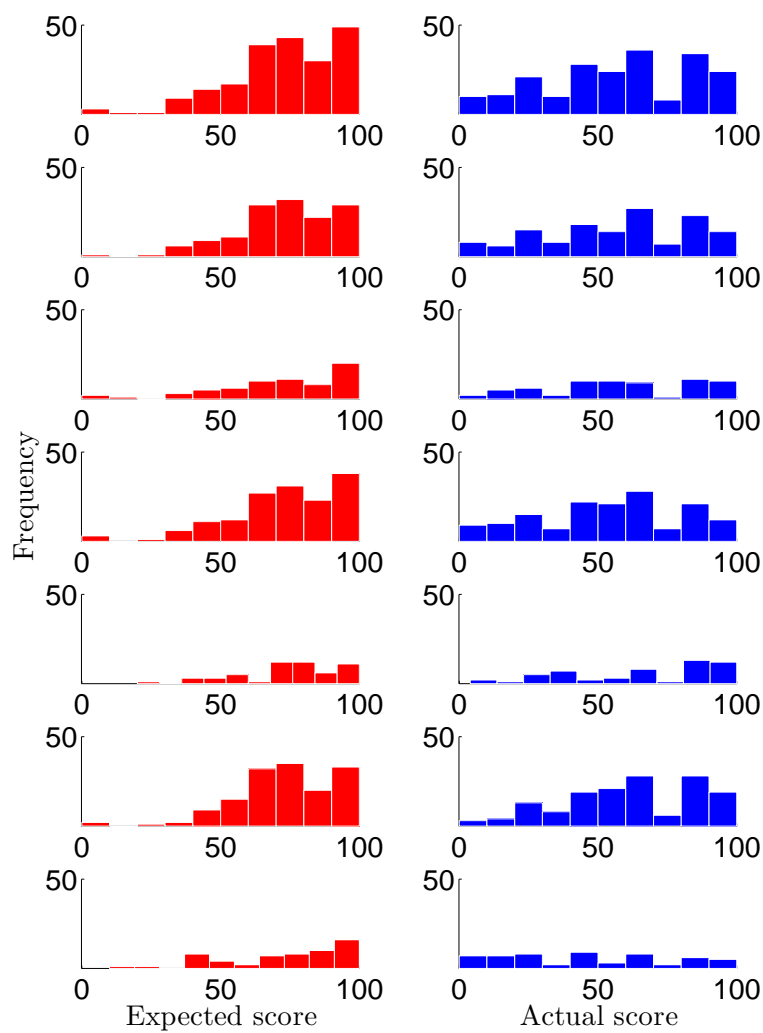}
  \caption{Histograms based on the collected sample data of the expected score (red histograms, left panels) and the actual AP (blue histograms, right panels) for all categories. From the top to the bottom, the histograms are arranged in the following order: all students, Foundation students, Undergraduate students, male students, female students, local students and international students, respectively.}
  \label{Histo2}
\end{center}
\end{figure}

The bar charts of data distribution for both the expected and the actual scores are presented as histograms in Figure~\ref{Histo1}. From the histogram for the expected scores, it can be observed that the data for the expectation on AP show a non-normal left-skewed type of distribution. This indicates that generally students tend to have a high expectation score regarding the test, irrespective of their actual academic ability. On the other hand, the histogram for the actual score distribution shows a plateau (multimodal) type of distribution where the data are distributed differently for each score interval. Furthermore, the bar chart distributions of the collected sample data for all categories for both the expected and the actual scores are presented as histograms in Figure~\ref{Histo2}. From the top to the bottom, the histograms are arranged in the following order: all students, Foundation students, Undergraduate students, male students, female students, local students and international students, respectively. The left panels (red histograms) present data distribution of the expected score and the right panels (blue histograms) present data distribution of the actual AP. From these histograms, it can also be observed that the students in each category generally tend to have a high expectation score regarding the test, while the actual score data are distributed randomly for each score interval.

It can be noticed that there is no significant correlation in the averages of the expected and the actual scores between the students enrolled in both the Foundation and Undergraduate programmes. A~similar situation is also observed by comparing the averages of the expected and the actual scores between the local and international students. However, by comparing the gender groups, even though there is no significant correlation in the average of the expected score, there is a significant difference ($p < 0.01$) between the actual score average of the male students and that of the female students.

The results of the responses to the seven statements from the distributed questionnaires are given in Table~\ref{percentage} and Table~\ref{average2}. Table~\ref{percentage} reports the number of students who responded agree and strongly agree to each questionnaire statement and its corresponding percentage. Table~\ref{average2} reports the average values of each questionnaire statement based on the Likert scale ranging from 1 to 4. Both tables list the information for all students, the students enrolled in the Foundation and Undergraduate programmes, and the male and female students as well as the local and international students. It is discovered that more than 70\% of the students (71.84\%) feel calm before the test and almost two third (65.53\%) feel confident enough to face the test. In addition, more than 40\% feel tense, more than a quarter feel panic (28.15\%) and are emotionally disturbed (26.70\%), more than half feel nervous (56.80\%) and 44.17\% of the population feel overexcited and at the same time feel worried.

By comparing the students based upon their academic levels, it is observed that there is no significant difference in the levels of test anxiety, except that the students enrolled in the Foundation programme were more nervous than the Undergraduate students, which is significant at 5\% level. By comparing between gender groups, it is perceived that generally female students have a higher level of test anxiety than male students. It is discovered that the male students are more confident and feel calmer than the female students, where these results are significant at 1\% and 5\% levels, respectively. The total average values for both gender groups are also significant at 5\% level. However, the male students are actually overconfident. There is at least 11\% difference between the two groups in terms of the actual AP, as can be seen in Table~\ref{average1}. It seems that the females are a better group of students, and thus it appears that they are better predictors of their success than their male counterparts. By comparing the local and international students, it is noticed that generally local Malaysian students have a higher level of test anxiety in mathematics than the international students. No significant difference between these two groups of students is observed.

\begin{table}[h]
{\footnotesize
\begin{center}
\caption{The number of students ($N$) who opted agree and strongly agree to each Likert questionnaire statement and the corresponding percentage~(\%) of all students, the students enrolled in the Foundation and Undergraduate programmes, the male and female students as well as the local and international students.} \label{percentage}
\begin{tabular}{lrcrcrcrcrcrcrc}
\hline \hline
 & \multicolumn{2}{c}{All} & \multicolumn{2}{c}{Foundation} & \multicolumn{2}{c}{Undergraduate} & \multicolumn{2}{c}{Male} & \multicolumn{2}{c}{Female} & \multicolumn{2}{c}{Local} & \multicolumn{2}{c}{International} \\
\cline{2-15}
Likert item & $N$ & \% & $N$ & \% & $N$ & \% & $N$ & \% & $N$ & \% & $N$ & \% & $N$ & \% \\
\hline \hline
Calm & 148 & 71.85 & 102 & 72.86  & 46  & 69.70  & 115 & 76.16  & 33 & 60.00 & 103 & 69.59 & 45 & 77.59 \\
Tense & 83  & 40.29 & 54  & 38.57  & 29  & 43.94  & 57  & 37.75  & 26 & 47.27 & 61  & 41.22 & 22 & 37.93 \\
Emotion & 55  & 26.70 & 35  & 25.00  & 20  & 30.30  & 37  & 24.50  & 18 & 32.73 & 38  & 25.68 & 17 & 29.31 \\
Nervous & 117 & 56.80 & 81  & 57.86  & 36  & 54.55  & 81  & 53.64  & 36 & 65.46 & 88  & 59.46 & 29 & 50.00 \\
Panic & 58  & 28.15 & 29  & 20.71  & 29  & 43.94  & 42  & 27.81  & 16 & 29.09 & 38  & 25.68 & 20 & 34.48 \\
Confident & 135 & 65.53 & 106 & 75.71  & 29  & 43.94  & 107 & 70.86  & 28 & 50.91 & 97  & 65.54 & 38 & 65.52 \\
Worry & 91  & 44.17 & 59  & 42.14  & 32  & 48.48  & 61  & 40.40  & 30 & 54.55 & 70  & 47.30 & 21 & 36.21 \\
\hline \hline
\end{tabular}
\end{center}
Notes: The Likert item `Calm' denotes `I feel calm', `Tense' denotes `I am tense', `Emotion' denotes `I am emotionally disturbed', `Nervous' denotes `I feel nervous', `Panic' denotes `I am panicky', `Confident' denotes `I am confident' and `Worry' denotes `I feel overexcited and worried'. \par}
\end{table}

The correlation values between each of the seven Likert questionnaire statements are listed in the correlation matrix shown in Table~\ref{matrix}. Since the matrix is symmetric, the upper triangular components have values that are the same as those of the lower triangular ones and are reflected with respect to the diagonal components. Only the lower triangular components are given in the table. It is observed that all Likert questionnaire statements have positive correlation values except for the first statement `I feel calm' and the sixth statement `I am confident'. These two statements have a positive correlation value between themselves. The absolute value of the correlation coefficients ranges from 0.202 to 0.621. The smallest correlation coefficient in the absolute value sense occurs between the statements `I am confident' and `I feel overexcited and worried' ($r = -0.202$). The effect size of this negative correlation value is small where there is a linear relationship between the two Likert items.\footnote{An `effect size' is a simple way of quantifying the size of the difference between two variables. In this article, the effect size of the correlation value from \citeauthor{Cohen88}~(\citeyear{Cohen88}) is adopted. The author proposed that $r = 0.1$, $r = 0.3$ and $r =0.5$ correspond to small, moderate/medium and large effect sizes, respectively. A moderate effect size represents an effect that is most likely visible to the naked eye of a careful observer. A small effect size is noticeably smaller than the moderate one, but not so small as to be trivial. A large effect size is grossly perceptible and has the same distance above the moderate one as the small one has below it (\citealp{Cohen92a,Cohen92b}).} The negative sign indicates that the students with more confidence would feel less overexcited and worried and vice versa. The largest correlation coefficient occurs between the statements `I feel nervous' and `I am panic' ($r = 0.621$). The effect size of this positive correlation value is considered large and suggests a fairly predictable relationship between the two statements. The positive sign indicates that the students who felt nervous would also panic and vice versa.

\begin{figure}[t]
\begin{center}
  \includegraphics[width=0.40\textwidth]{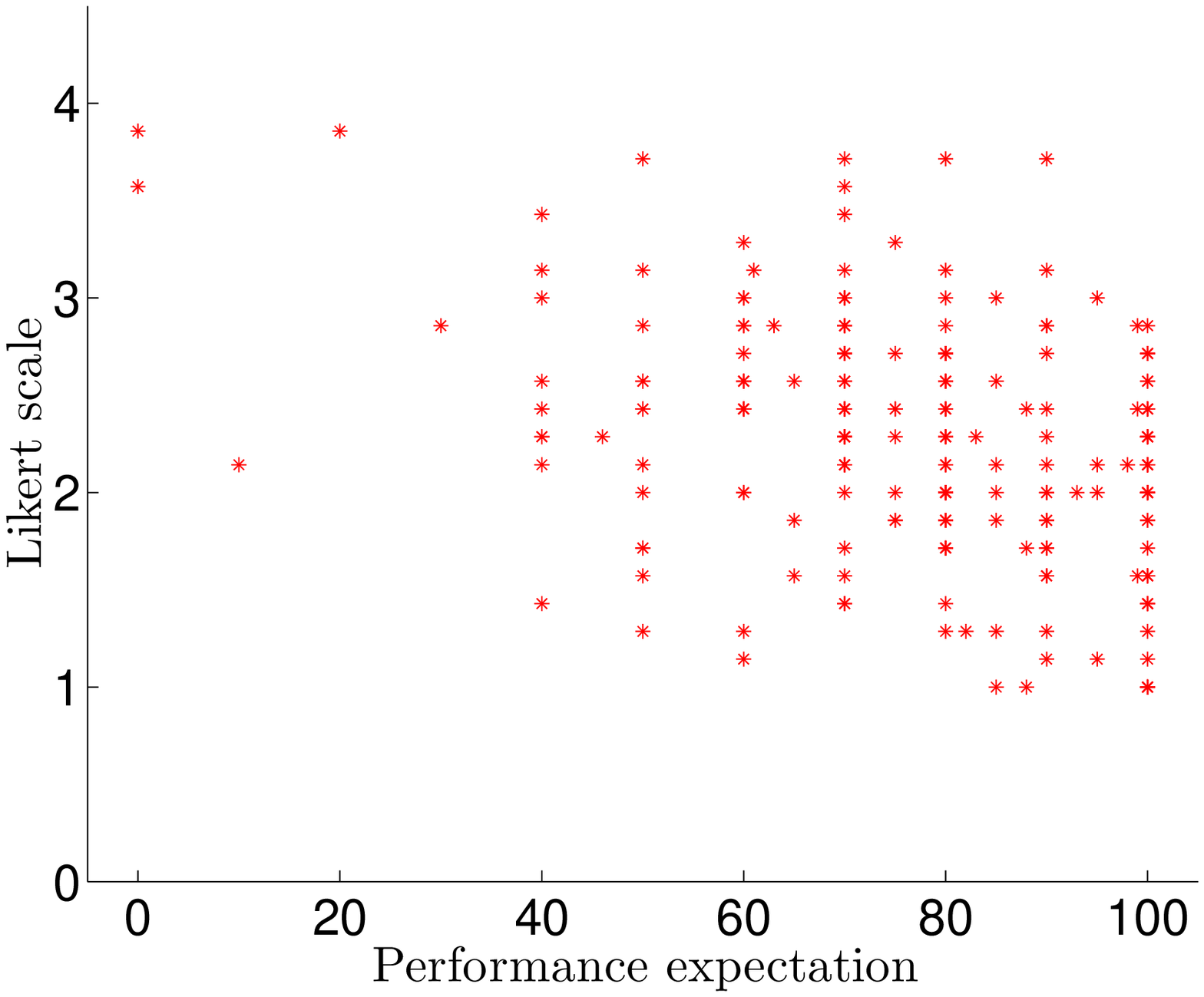} \hspace{0.5cm}
  \includegraphics[width=0.45\textwidth]{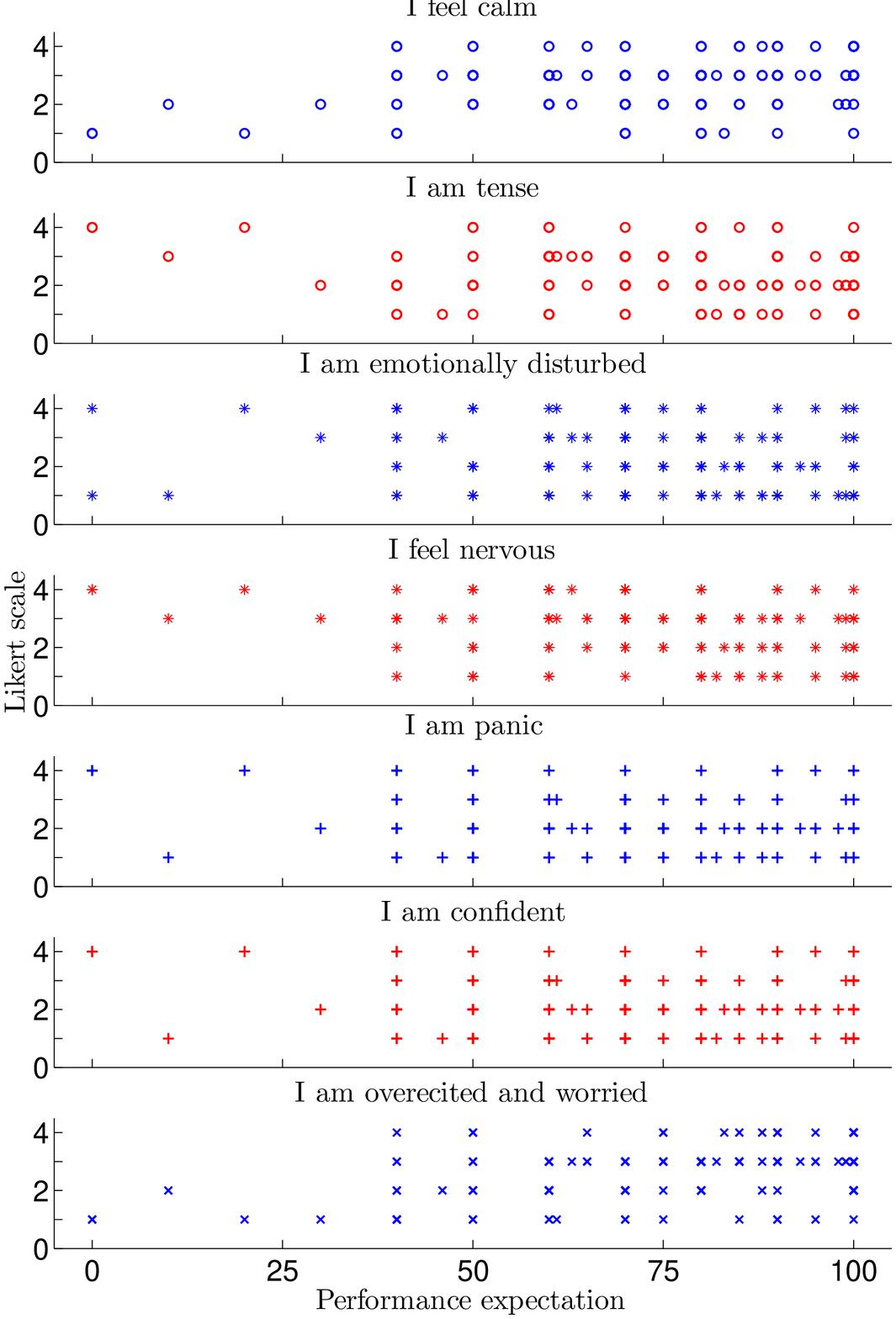} \\
  \caption{Scatter plots of the level of test anxiety in mathematics versus the expectation on AP. In the left panel (red scatter), the vertical axis denotes the average score of each Likert questionnaire item based upon the collected sample data for each student and the horizontal axis is the students' expected score during a particular test. The right panel presents the scatter plots for all seven Likert items versus the expected score.}
  \label{NewScatter0}
\end{center}
\end{figure}

\begin{figure}[h]
\begin{center}
  \includegraphics[width=0.40\textwidth]{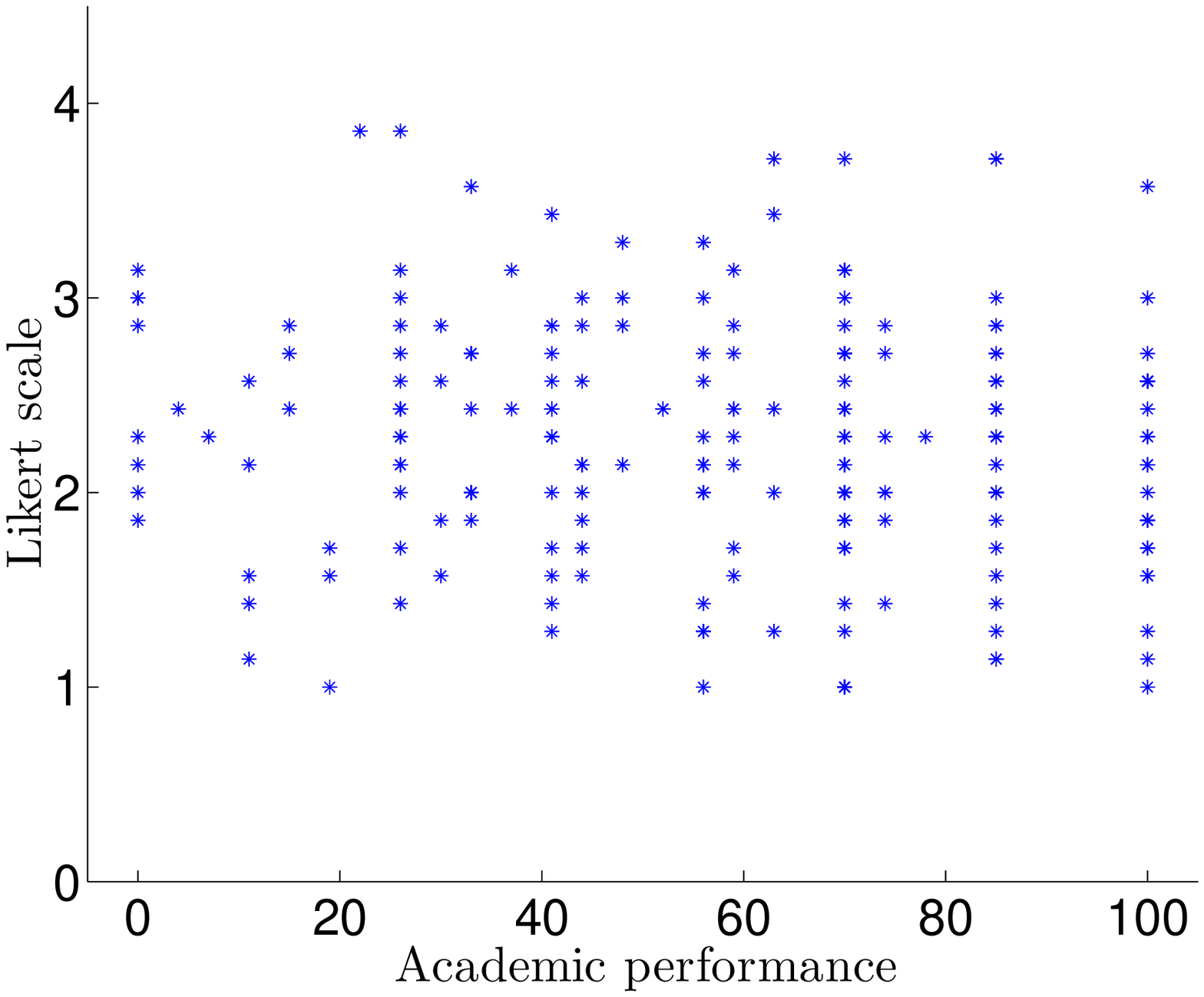} \hspace{0.5cm}
  \includegraphics[width=0.45\textwidth]{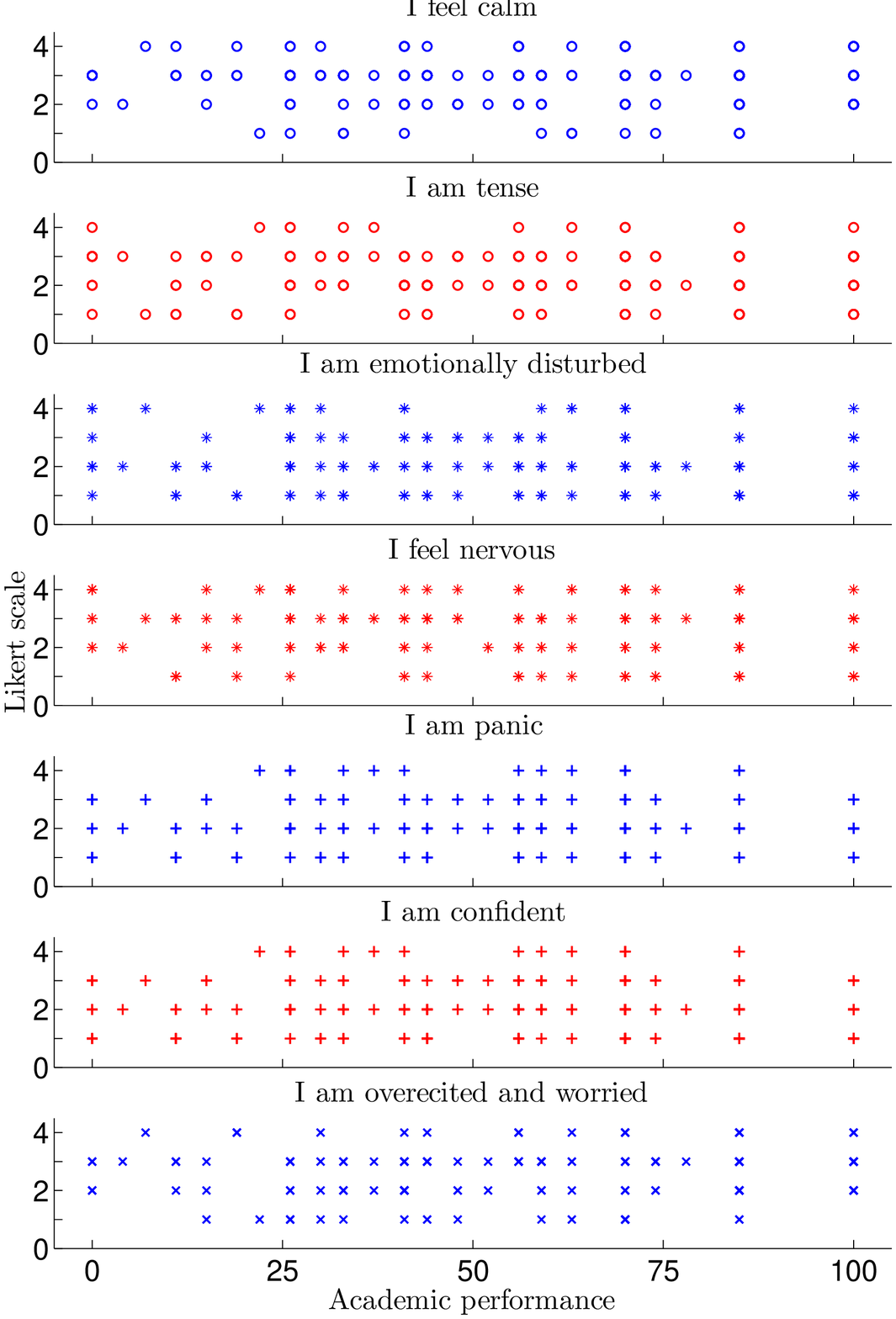} \\
  \caption{Scatter plots of the level of test anxiety in mathematics versus the actual AP. In the left panel (blue scatter), the vertical axis denotes the average score of each Likert questionnaire item based upon the collected sample data for each student and the horizontal axis denotes the actual score that the students obtained for the particular test. In the right panel, the scatter plots for all seven Likert items versus the actual AP are presented.}
  \label{NewScatter1}
\end{center}
\end{figure}

\begin{figure}[h]
\begin{center}
  \includegraphics[width=0.45\textwidth]{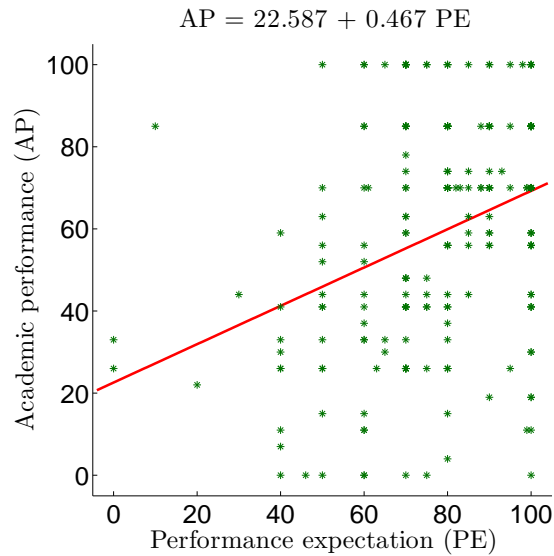} \\
  \caption{A scatter plot (green scatter) of the expected score (horizontal axis) and the actual AP (vertical axis). A simple linear regression to the data gives a linear relationship (red line) of AP $ = 22.587 + 0.467$ PE, where AP and PE denote the students' actual score and expectation on AP, respectively.}
  \label{NewScatter2}
\end{center}
\end{figure}

The visualizations of the collected sample data are presented in the scatter plots shown in Figures~\ref{NewScatter0}--\ref{NewScatter2}. The scatter plots of the collected data for the average score of each Likert item for all students versus the expectation on AP are presented in Figure~\ref{NewScatter0}. In the left panel (red scatter), it can be observed that many students have a high expectation on AP irrespective of the level of test anxiety. In the right panel, the scatter plots for all seven Likert statement items versus the expected score distribution are displayed. A similar tendency is also observed that the majority of students hope to achieve high scores during the tests. The scatter plots of the level of test anxiety versus the actual AP in mathematics are shown in Figure~\ref{NewScatter1}. In the left panel (blue scatter), the actual scores are more or less evenly distributed for varied numbers of the average level of test anxiety, except that there are no scores between 85 and 100. In the right panel, the scatter plots for all seven Likert statement items versus the actual AP data are displayed. Except for the scores between 85 and 100, the score data are distributed throughout the possible score range, irrespective of the students' level of test anxiety in mathematics.
\begin{table}[h]
{\footnotesize \renewcommand{\tabcolsep}{8pt}
  \begin{center}
  \begin{minipage}{16cm}
  \caption{The average values of each Likert questionnaire statement item for all students, the students enrolled in the Foundation and Undergraduate programmes, the male and female students as well as the local and international students.} \label{average2}
    \begin{tabular}{lccccccc}
      \hline \hline
      Statement item & All & Foundation & Undergraduate & Male & Female & Local & International \\
      \hline \hline
      I feel calm$^{a}$              & 2.22 & 2.23 & 2.20 & 2.15 & 2.40$^{c}
      \!\!\!$ & 2.24 & 2.16 \\
      I am tense                     & 2.31 & 2.31 & 2.29 & 2.26 & 2.44 & 2.33 & 2.56 \\
      I am emotionally disturbed     & 2.04 & 2.12 & 1.82 & 1.98 & 2.22 & 2.04 & 2.05 \\
      I feel nervous                 & 2.56 & 2.64 & 2.32$^{c}\!\!$ & 2.50 &
      2.73 & 2.59 & 2.48 \\
      I am panic                     & 2.09 & 2.06 & 2.16 & 2.05 & 2.18 & 2.07 & 2.14 \\
      I am confident$^{a}$           & 2.33 & 2.37 & 2.20 & 2.25 & 2.55$^{b}
      \!\!\!$ & 2.38 & 2.19 \\
      I feel overexcited and worried & 2.36 & 2.40 & 2.25 & 2.32 & 2.47 & 2.41 & 2.24 \\
      Total average                  & 2.27 & 2.30 & 2.18 & 2.22$^{c}\!\!$ &
      2.43$^{c}\!\!$ & 2.30 & 2.22 \\ \hline \hline
    \end{tabular}\\

    Notes: The possible minimum and maximum values are 1 and 4, respectively.
    \footnotetext[1]{The corresponding scores are reversed to calculate the total average.}
    \footnotetext[2]{Significant at 0.01 level.}
    \footnotetext[3]{Significant at 0.05 level.}
  \end{minipage}
\end{center}
}
\end{table}
\begin{table}[h]
{\footnotesize \renewcommand{\tabcolsep}{13pt}
\begin{center}
\caption{A table of intercorrelation values (correlation matrix) among the seven Likert items.}\label{matrix}
  \begin{tabular}{lrrrrrrr}
    \hline \hline
    Likert item  & Calm & Tense & Emotion & Nervous & Panic & Confident & Worry \\ \hline \hline
    Calm & 1{\color{white}.000} &   &   &   &   &   &  \\
    Tense & $-0.528$ & 1{\color{white}.000} &   &   &   &   &  \\
    Emotion & $-0.315$ & $0.330$  & 1{\color{white}.000} &   &   &   &  \\
    Nervous & $-0.533$ & $0.537$  & $0.476$  & 1{\color{white}.000} &   &   &  \\
    Panic & $-0.534$ & $0.527$  & $0.549$  & $0.621$  & 1{\color{white}.000} &   &  \\
    Confident & $0.423$  & $-0.272$ & $-0.351$ & $-0.459$ & $-0.469$ & 1{\color{white}.000} &  \\
    Worry & $-0.358$ & $0.344$  & $0.446$  & $0.359$  & $0.414$  & $-0.202$ & 1{\color{white}.000} \\
    \hline \hline
  \end{tabular}
\end{center}
Notes: Since the correlation matrix is symmetric, only the lower triangular components are listed. The upper triangular components have symmetric values. The Likert items `Calm' denotes `I feel calm', `Tense' denotes `I am tense', `Emotion' denotes `I am emotionally disturbed', `Nervous' denotes `I feel nervous', `Panic' denotes `I am panicky', `Confident' denotes `I am confident' and `Worry' denotes `I feel overexcited and worried'. \par}
\end{table}

As the third side of the correlation triangle, a scatter plot between the students' expectation on the test and the actual AP is shown in Figure~\ref{NewScatter2}.
It can be seen that the collected sample data for the two types of scores tend to scatter towards the right part of the figure, which indicates that the students generally expect a high score than what they could actually achieve in the test. A positive correlation between the two variables is also clearly noticed where the students with a higher expectation on the test scores will also exhibit higher actual AP. Additionally, a simple linear regression to the collected sample data shows a linear relationship between the expected score (PE) and the actual AP, given as AP $= 22.587 + 0.467$ PE. Assuming that expected score (PE) is a predictor variable and the actual AP is a response variable, the AP-intercept of 22.587 indicates the value of an actual score if a student expects zero score during the test. The slope of 0.467 indicates the amount that the actual score will increase as the expected score increases by one point. According to this linear relationship, an average score of 76.44 for PE would predict an average of 58.285 for the actual AP. This prediction is almost accurate as the average score of the actual AP from the collected data is 58.26. It shows that, on average, the students overestimate their scores by about 18\%.

The correlations between the students' expectation score and the level of test anxiety and those between the level of test anxiety and the AP for all students, students in different study programmes, male and female students and students from different nationality backgrounds (local versus international) are investigated Subsections~\ref{all}--\ref{cultural}.

\subsection{Entire group analysis} \label{all}

In order to test Hypothesis~\ref{hypothesis1} (there is a significantly negative correlation between the students' expectation on AP in the mathematics test and the level of test anxiety in mathematics during the test period), the correlation value between the students' expectation score and their level of mathematics-test anxiety is calculated. The correlation values $r$ and the corresponding levels of significant $p$-values for all the students are listed in Table~\ref{correlation2}. It is noted that there is a significant ($p < 0.01$) negative relationship ($r = -0.366$) between the expectation score and the level of test anxiety in mathematics during the test period. This means that the students with a higher expectation on their AP in the test will experience a lower level of test anxiety in mathematics. The effect size of the correlation value is moderate, the relationship is not really strong and yet there is a definite linear relationship between the level of mathematics-test anxiety and the expected score. Since there is a significant correlation, this result supports Hypothesis~\ref{hypothesis1}.

In order to test Hypothesis~\ref{hypothesis2} (there is a significantly negative correlation between the students' level of test anxiety in mathematics during the test period and their AP in mathematics), the correlation between the students' level of mathematics-test anxiety and their actual scores obtained in the test is calculated. It is also discovered that there exists a negative relationship ($r = -0.099$) between the level of mathematics-test anxiety and the actual score, even though the result is not significant ($p = 0.157$). This means that the students with a higher level of mathematics-test anxiety will perform worse in the test than those with a lower level of mathematics-test anxiety. The effect size of this correlation value is very small and can be interpreted as a very small linear relationship between the level of test anxiety and the AP in a negative manner. Since there is no significant correlation, this finding does not support Hypothesis~\ref{hypothesis2}.

In order to test Hypothesis~\ref{hypothesis3} (there is a significantly positive correlation between the students' expectation on AP in the mathematics test and their actual AP in mathematics), the correlation value between students' level of mathematics-test anxiety and their actual test scores is calculated. Different from the previous two cases, there exists a positive relationship ($r = 0.333$) between the expectation score and the actual score. This result is also significant at 1\% level ($p = 0.000$). This indicates that the students with a high expectation will exhibit better AP in mathematic tests. The effect size of this correlation value is considered as moderate in psychological research. Although the relationship is not that strong, there is definitely a linear relationship between the expectation score and the actual score variables. Since the correlation is significant, this finding supports Hypothesis~\ref{hypothesis3}.

\subsection{Study programme comparison} \label{programme}

The students who enrol in the Foundation programme are those who have completed their secondary education with the ordinary level (O~level) qualification of the General Certificate of Education (GCE) or equivalent. On the other hand, the students who enrol in the Undergraduate programme are those who have completed their high school education with the advanced level (A~level) qualification of the GCE or equivalent.\footnote{The GCE is an academic qualification that examination boards confer to students in the UK and a few of the commonwealth countries, including Hong Kong, Malaysia, Singapore, Sri Lanka and Pakistan. Traditionally, it comprises two levels: the O~Level and the A~Level. More recently, examination boards have also offered an intermediate third GCE level, the Advanced Subsidiary Level, replacing the earlier Advanced Supplementary level.} In fact, the majority of students enrolled in the Undergraduate programme at the UNMC are also formerly enrolled in the Foundation programme at the same campus. Only a small percentage of students have completed a similar Foundation programme somewhere else in either Malaysian or overseas tertiary institutes and with direct entry with the A~level qualification or equivalent.

In this subsection, Hypothesis~\ref{hypothesis4} (there is a significant difference in AP between the students enrolled in the pre-undergraduate (Foundation) programme and the ones enrolled in the Undergraduate programme when dealing with test anxiety in mathematics) is examined. It is observed from Table~\ref{percentage} that the percentage of students enrolled in the Foundation programme who showed some level of mathematics-test anxiety ranges from around one-fifth for the fifth Likert item `I am panicky' (20.71\%) to more than half for the fourth Likert item `I am nervous' (57.86\%).  On the other hand, the percentage of students enrolled in the Undergraduate programme who showed some level of mathematics-test anxiety ranges from around one-third for the first (reversed) and the third Likert items `I do not feel calm' and `I am emotionally disturbed', respectively (30.30\%) to more than half for the sixth Likert item (reversed) `I am not confident' (56.06\%). So far, no clear conclusion has been drawn from this percentage figures.

The average values of the mathematics-test anxiety level for each Likert questionnaire statement for both the Foundation and Undergraduate students are calculated and listed in Table~\ref{average2}. It is discovered that for the students enrolled in the Foundation programme, the average values for each Likert item range from 2.06 to 2.64 for the fifth and the fourth Likert items, respectively. The average values for the students enrolled in the Undergraduate programme range from 1.82 to 2.32 for the third and the fourth Likert items, respectively. The average value of 2.32 is statistically significant at 5\% level. The total average values of the level of test anxiety for the Foundation and Undergraduate students are 2.30 and 2.18, respectively. These indicate that the students enrolled in the Foundation programme have a higher level of test anxiety in mathematics in comparison with those enrolled in the Undergraduate programme.

\begin{table}[h]
{\footnotesize \renewcommand{\tabcolsep}{9pt}
\begin{center}
  \begin{minipage}{16cm}
  \caption{The values of correlation $r$ and the corresponding significant level $p$ for all students and the students enrolled in the Foundation and Undergraduate programmes.}
  \label{correlation2}
\begin{tabular}{lccc}
\hline \hline
  & \multicolumn{3}{c}{\hspace*{0cm} Programme enrolled} \\
\cline{2-4}
Correlation of data variables & All & Foundation & Undergraduate \\
\hline
Expectation score versus level  & $r = -0.366^{d}$ ($p = 0.000$) & $r = -0.411^{d}$ ($p = 0.000$)  & $r = -0.210$ $(p = 0.143)$ \\
of mathematics-test anxiety     &       &        &  \\
Test anxiety level versus  & $r = -0.099$ ($p = 0.157$) & $r = -0.103$ ($p = 0.200$) & $r = -0.029$ ($p = 0.842$) \\
actual score               &        &       &      \\
Expectation score versus   & $r = 0.332^{d}$ ($p = 0.000$)  & $r = 0.333^{d}$ ($p = 0.000$) & $r = 0.337^{e}$ $(p = 0.017)$ \\
actual score               &        &       &   \\
\hline \hline
\end{tabular}
\footnotetext[4]{Significant at the 0.01 level.}
\footnotetext[5]{Significant at the 0.05 level.}
\end{minipage}
\end{center}
}
\end{table}

In order to test the null hypothesis as to whether the means between the two groups of students are equal against the non-directional alternative hypothesis that the two means differ, the two-tailed Student $t$-test for two independent samples is implemented. However, the $t$-test is valid under the assumption that the distribution of data in the underlying population from which each sample is  derived is normal. For the sample groups that do not appear to be normally distributed, an alternative non-parametric statistical hypothesis test of Wilcoxon rank-sum $W$-test or an identical Mann-Whitney $U$-test should be implemented. These two versions of the test are identical, and they were independently developed by \citeauthor{Wilcoxon45}~(\citeyear{Wilcoxon45}) and \citeauthor{Mann47}~(\citeyear{Mann47}). Both differ in terms of the constant factor of the mean (\citealp{Hollander73}). The results of the Student $t$-test, Wilcoxon rank-sum $W$-test and Mann-Whitney $U$-test and their corresponding $p$-values are presented in Table~\ref{tWU}.

Taking an assumption that the sample groups are normally distributed, the Student $t$-test is implemented. It is observed that the $t$ value for the two-tailed $t$-test is 0.736. Since this value is less than the critical values $t_{.05} = 1.97$ and $t_{.01} = 2.60$ with the degree of freedom $205$, the null hypothesis that the means are equal is accepted and thus there is no difference between the two groups. Furthermore, since the $p$-value is 0.462, there is no significant difference between the students enrolled in the Foundation programme and those enrolled in the Undergraduate programme. Since it is not known that the data are normally distributed, either the Wilcoxon rank-sum $W$-test or the Mann-Whitney $U$-test can be implemented to test the null hypothesis whether both groups of students have the same median against the non-directional alternative hypothesis that the two groups have different median values. According to the Table of the Normal Distribution, the two-tailed 0.05 and 0.01 critical values are $z_{.05} = 1.96$ and $z_{.01} = 2.58$, respectively. Since the absolute values of both tests exceed these critical values ($z = -18.687$, $p = 0.000$ for $W$-test and $z = 6.108$, $p=0.000$ for $U$-test), the null hypothesis is rejected at both level of significance. Thus, the median of the two groups of students are not identical and the difference between them is significant.

Nonetheless, the effect of test anxiety on the AP in mathematics may differ from one group to the other. The analysis goes further to determine the correlation between the expectation score and the level of mathematics-test anxiety as well as the correlation between the level of mathematics-test anxiety and the actual score among the students enrolled in the Foundation and Undergraduate programmes. This finding is then used to examine whether there is a significant correlation between the level of mathematics-test anxiety and the AP in the context of different levels of study. The correlation values $r$ and the corresponding levels of significance $p$-values for the students enrolled in the Foundation and Undergraduate programmes are given in Table~\ref{correlation2}. It is discovered that the level of mathematics-test anxiety for the Foundation students is significantly related to their expectations in a negative manner, with correlation value of $r = -0.411$ at a significant level of 1\% ($p = 0.000$). The effect size of this correlation value is moderate and there is a negative linear relationship between the level of test anxiety and the expected score.

\begin{table}[h]
{\footnotesize \renewcommand{\tabcolsep}{13pt}
\begin{center}
\caption{The values of the Student $t$-test, Wilcoxon rank-sum $W$-test and Mann-Whitney $U$-test together with their corresponding significant level $p$-values.} \label{tWU}
\begin{tabular}{lcccccc}
\hline \hline
    &
\multicolumn{6}{c}{\hspace*{0cm} Statistical tests} \\
\cline{2-7}
Student categories & $t$-Test & $p$-Value & $W$-Test & $p$-Value &$U$-test & $p$-Value \\
\hline
Foundation and Undergraduate  & 0.736 & 0.462 & $-18.687$ & 0.000 & 6.108 & 0.000 \\
Male and female               & 0.032 & 0.975 & $-24.043$ & 0.000 & 6.366 & 0.000 \\
Local and international       & 0.412 & 0.681 & $-25.036$ & 0.000 & 3.702 & 0.000 \\
\hline \hline
\end{tabular}
\end{center}
Notes: The ties correction has been incorporated into the latter two tests.
}
\end{table}

Furthermore, a negative correlation $(r = -0.103)$ is also observed between the level of mathematics-test anxiety and the obtained actual score, but the difference is not significant ($p = 0.200 > 0.05$). The effect size of this correlation value is very small. There is a very small relationship in a negative manner between the level of test anxiety and the AP. The third side of the correlation triangle shows that there is a positive correlation ($r = 0.333$) between the expectation score and the actual score with a significant $p$-value at 1\% level ($p = 0.000$). The effect size of this correlation value is considered moderate, indicating that there is a linear relationship between the two variables. Since the number of students enrolled in the Foundation programme is more than two times the number of those enrolled in the Undergraduate programme ($140:66$), it is not surprising that the obtained result would be similar to the entire group comparison.

Moreover, any significant association between the level of test anxiety and the expected as well as the actual scores among the students enrolled in the Undergraduate programme is not noticed. Negative correlation values are observed for both cases. The correlation value for the level of mathematics-test anxiety and the expected score is $r = -0.210$ and the difference is not significant ($p =0.144 > 0.05$). The effect size is small, which means that a high expectation score would predict a low level of test anxiety but not with a great deal of accuracy. Thus, the expectation has little or no effect on the level of test anxiety during the time of the midterm examination. Likewise, the state of test anxiety has no significant impact on the AP in mathematics. The correlation value for the level of mathematics-test anxiety and the AP is $r = -0.029$ and the result is not significant either ($p = 0.842 > 0.05$). Since the absolute value of $r$ is close to zero, the effect size is very small. This informs that almost no relationship between the two variables. Since there is no significant correlation, this result does not support Hypothesis~\ref{hypothesis4}.

The third side of the correlation triangle again shows a positive correlation value ($r = 0.337$) between the expected score and the obtained actual score for the students enrolled in the Undergraduate programme. The $p$-value was found to be $p = 0.017$ and is significant at 5\% level. The effect size of the correlation value is moderate, and this indicates that there is a linear relationship between the two variables of the expected score and the actual score.

\subsection{Gender difference comparison} \label{gender}

To check Hypothesis~\ref{hypothesis5} (there is a significant difference in AP between male and female students when dealing with test anxiety in mathematics), the percentages of each Likert item and the average values of the level of mathematics-test anxiety for both male and female students are also calculated. It is noticed from Table~\ref{percentage} that the percentage of male students who showed some level of mathematics-test anxiety ranges from almost a quarter for the first (reversed) Likert item `I am not calm' (23.84\%) to more than half for the fourth Likert item `I am nervous' (53.64\%). The percentage of female students ranges from almost one-third for the third Likert item `I am emotionally disturbed' (32.73\%) to almost two-thirds for the fourth Likert item `I am nervous' (65.46\%).

The average values for each Likert item for the male and female students can be seen in  Table~\ref{average2}. The average values for the male students range from 1.98 (the third Likert item) to 2.50 (the fourth Likert item) and those for the female students range from 2.18 (the fifth Likert item) to 2.55 (the sixth Likert item). The average value of 2.55 for the sixth Likert item is statistically significant at 1\% level. The total average values of the level of mathematics-test anxiety for the male and female students are 2.22 and 2.43, respectively. From these figures, it seems that trhe female students have a higher level of mathematics-test anxiety than the male students.

Furthermore, the two tailed Student $t$-test result gives the $t$-value of 0.032. This indicates that there is no difference between the groups of male and female students. Since the $p$-value is 0.975, the $t$-test result shows that there is no significant difference between these two groups. However, both the Wilcoxon rank-sum $W$-test ($z = -24.043$, $p = 0.000$) and the Mann-Whitney $U$-test ($z = 6.366$, $p=0.000$) give us the information that there is a significant difference at 1\% level between the two groups of male and female students.
\begin{table}[h]
{\footnotesize \renewcommand{\tabcolsep}{10pt}
\begin{minipage}{16cm}
\caption{The values of correlation $r$ and the corresponding significant level $p$-values in the context of gender difference.} \label{correlation3}
\begin{tabular}{lcc}
\hline \hline
 &
\multicolumn{2}{c}{\hspace*{0cm} Gender} \\
\cline{2-3}
Correlation of data variables & Male & Female \\
\hline
Expectation score versus level of mathematics-test anxiety  & $r = -0.390^{f}$ ($p = 0.000$)  & $r = -0.283^{g}$ ($p = 0.036$) \\
Test anxiety level versus actual score                      & $r = -0.169^{f}$ ($p = 0.038$)  & $r = -0.013$     ($p = 0.925$) \\
Expectation score versus actual score                       & $r =  0.325^{f}$ ($p = 0.000$)  & $r =  0.399^{f}$ ($p = 0.000$) \\
\hline \hline
\end{tabular}
\footnotetext[6]{Significant at the 0.01 level.}
\footnotetext[7]{Significant at the 0.05 level.}
\end{minipage}
}
\end{table}

Nevertheless, the interest is on whether the level of test anxiety on both gender groups exerts the same effect on their AP in mathematics. The analysis is carried out further to determine the correlation between the expectation score and the level of mathematics-test anxiety as well as the correlation between the level of mathematics-test anxiety and the actual score between the male and female students. The values of correlation $r$ and the corresponding significant level $p$-values of the collected data in the context of gender difference are given in Table~\ref{correlation3}.

It is noticed that the level of mathematics-test anxiety is correlated significantly with the expectation score in a negative manner for both gender groups. The correlation values for the male and female students are $r = -0.390$ and $r = -0.283$ with significant at 1\% and 5\% levels, respectively. These indicate that both gender groups with a higher expectation score would experience a lower level of mathematics-test anxiety. The effect sizes of both correlation values can be considered as moderate. However, since the absolute value of the correlation value for the male students is larger than that of the correlation value for the female students, a stronger negative relationship between the expectation score and the level of mathematics-test anxiety is more evident among the male students than among the female students.

Moreover, negative correlations are also observed between the level of mathematics-test anxiety and the actual AP, with the correlation value for the male students $r = -0.169$. Thus, the male students who possess a higher level of mathematics-test anxiety would exhibit lower AP. The effect size of the correlation value for the male students is small and hence a small negative relationship would be observed between the level of test anxiety and the AP. The effect size of the correlation value for the female students is extremely small and this implies that almost no relationship exists between the two variables. It is noticed that only male students have a significant correlation at 5\% level ($p = 0.038$) and the result for the female students is not significant ($p = 0.925 > 0.05$). Since there is no significant correlation for the female students, this finding does not support Hypothesis~\ref{hypothesis5}.

The third side of the correlation triangle between the expectation score and the actual AP shows similar positive correlation values for both the male and female students. This is interesting since the number of male students is around 2.75 times the number of female students (151:55). The correlation values for the male and female students are $r = 0.325$ ($p = 0.000$) and $r = 0.399$ ($p = 0.000$), respectively. Both these correlation values are significant at 1\% level. This result indicates that both the male and female students who have a higher expectation score would also exhibit better AP. The effect sizes for both correlation values are moderate and give an indication that there are linear relationships between the expected score and the actual score for both the male and female students.

\subsection{Local versus international comparison} \label{cultural}

In this subsection, Hypothesis~\ref{hypothesis6} (there is a significant difference in AP between local Malaysian and international students when dealing with test anxiety in mathematics) is tested. The percentages of each Likert item and the averages values of the level of test anxiety for both local Malaysian and international students are calculated and are presented in Tables~\ref{percentage} and~\ref{average2}. It can be observed from Table~\ref{percentage} that the percentage of local Malaysian students who showed some level of mathematics-test anxiety ranges from around a quarter of the population (25.68\%) for the third (`I am emotionally disturbed') and the fifth (`I~am panic') Likert items to almost 60\% for the fourth Likert item `I am nervous' (59.46\%). The percentage of international students ranges from more than one-fifth (22.41\%) for the (reversed) first Likert item `I~am not calm' to half of the population (50\%) for the fourth Likert item.

It can be discovered from Table~\ref{average2} that the average values of each Likert item for local Malaysian students range from 2.04 (the third Likert item) to 2.59 (the fourth Likert item). The average values for the international students range from 2.05 (the third Likert item) to 2.48 (the fourth Likert item). None of these average values is statistically significant. The total average values for the local and international students are 2.30 and 2.22, respectively. From this finding, it seems that the local Malaysian students have a higher level of mathematics-test anxiety than  their overseas counterparts. This is a result that is rather different from the findings reported by \citeauthor{Burns91}~(\citeyear{Burns91}).

The Student $t$-test, Wilcoxon rank-sum $W$-test and Mann-Whitney $U$-test results for the local and international students are listed in Table~\ref{tWU}. It can be observed from Table~\ref{tWU} that the two-tailed $t$-test result for these two groups of students yields the $t$-value of 0.412. This indicates that there is no difference between the two groups. Furthermore, since the $p$-value is 0.681, the difference between the local and the international students is not significant. In addition to the $t$-test, the alternative non-parametric statistical tests for non-normal sample groups are also implemented. It can be found from Table~\ref{tWU} that the $W$-test gives $z = -25.036$ with $p$-value = 0.000 and the $U$-test gives $z = 3.702$ with $p$-value = 0.000. These results indicate that there is a significant difference at 1\% level between the two sample groups of local and international students.

Moreover, whether the level of test anxiety on these two different groups exerts a similar influence on their AP is examined. The investigation goes further by finding the correlation values between the expectation score and the level of mathematics-test anxiety as well as those between the level of test anxiety and the actual score obtained. The results, including the corresponding significant level $p$-values of the data, are presented in Table~\ref{correlation4}. It is noticed that negative correlation values are obtained between the expectation score and the level of mathematics-test anxiety for both groups of local Malaysian and international students. The corresponding correlation values for each group are $r = -0.117$ and $r = -0.209$, respectively. There is no significant correlation for both groups since the $p$-values for the local and international students are 0.157 and 0.115, respectively ($p > 0.05$). These figures show that both groups with a higher expectation on their AP would experience a lower level of mathematics-test anxiety. The effect sizes of the correlation values for local Malaysian and international students are small. Negative linear relationships can be observed between the expected score and the level of mathematics-test anxiety for both groups of students.
\begin{table}[h]
{\footnotesize \renewcommand{\tabcolsep}{10pt}
\begin{minipage}{16cm}
\caption{The values of correlation $r$ and the corresponding significant level $p$-values in the context of local versus international students.} \label{correlation4}
\begin{tabular}{lcc}
\hline \hline
 &
\multicolumn{2}{c}{\hspace*{0cm} Local versus international} \\
\cline{2-3}
Correlation of data variables & Local Malaysians & International students \\
\hline
Expectation score versus level of mathematics-test anxiety  & $r = -0.117$    ($p = 0.157$) & $r = -0.209$     ($p = 0.115$) \\
Test anxiety level versus actual score                      & $r = 0.037$     ($p = 0.655$) & $r = -0.396^{h}$ ($p = 0.002$) \\
Expected score versus actual score                          & $r = 0.326^{h}$ ($p = 0.000$) & $r =  0.344^{h}$ ($p = 0.008$) \\
\hline \hline
\end{tabular}
\footnotetext[8]{Significant at the 0.01 level.}
\end{minipage}
}
\end{table}

It is interesting to mention that a positive correlation value is obtained between the level of mathematics-test anxiety and the actual score for the local students ($r = 0.037$) while a negative correlation value is obtained for the international students ($r = -0.396$). The correlation is not significant for the former group of students ($p = 0.655 > 0.05$), but it is significant at 1\% level for the latter one ($p = 0.002$). Since the correlation value for the local students is positive, a higher level of mathematics-test anxiety would basically imply better AP. However, since the effect size of the correlation value for them is very small, it means that there is almost no relationship between the level of mathematics-test anxiety and the AP among the local students. On the other hand, for the international students, the correlation value is negative and its effect size is moderate. This means that a higher level of mathematics-test anxiety would imply lower AP and a linear relationship would be observed between these two variables. Since there is no significant correlation for the local Malaysian students, this finding does not support Hypothesis~\ref{hypothesis6}.

The third side of the correlation triangle between the expected score and the actual score demonstrated no surprise outcome since the results are similar to those obtained in the former cases. Positive correlation values are obtained for both sample groups of local and international students with $r = 0.326$ and $r = 0.344$, respectively. The ratio of local students to international students is around 2.5 (148:58), while the correlation values are similar up to the first digit of decimal. Both correlation values are significant at the 1\% level since the $p$-values for the local and international students are 0.000 and 0.008, respectively. This finding implies that both groups of local and international students with a higher expectation score would also exhibit better AP in mathematics. The effect sizes for both correlation values are considered moderate, and although the relationship is not strong, there is definitely a linear relationship between the expected score and the actual score for both local Malaysian and international students.

\section{Discussion} \label{conclusion}

The main purpose of this research is to investigate the level of test anxiety in connection to AP for several mathematics modules among early undergraduate students at the UNMC. This section presents some limitations of the study, draws some conclusions and gives future implications of the study.

\subsection{Limitations}

There exists a number of limitations in this research that may limit the generalisation of the findings to different populations. First, although the main target respondents were the early undergraduate students at the UNMC, the collected sample data were gathered from Engineering students who enrolled in the Undergraduate programme as well as from the students who were still enrolled in the Foundation programme in Engineering. The data regarding the students who were enrolled in other faculties were not collected. Even though these students also took several mathematics modules, they only did so at the Foundation level, but did not do it after they enrolled in the Undergraduate programme since some mathematics modules were no longer part of their curriculum. On the other hand, the students who enrolled in the Foundation program in Engineering would generally continue to take a number of mathematics modules after they progress into the Undergraduate level.

Second, the information regarding age, cumulative average points (similar to the GPA in the US educational system), socioeconomic status and ethnicity is not collected, as considered in \citeauthor{Chapel05}~(\citeyear{Chapel05}). The main reason for this is brevity; that is, the students could fill in the questionnaire in a relatively short period of time. This aspect is very essential since they were already under pressure to face the test. Third, the questionnaires adapted from TAI (\citealp{Spielberger80}) do not present a clear distinction between the measurements of the cognitive concerns and those of the emotional responses of the students associated with evaluation stress. Fourth, the composition of polled students is somewhat imbalance. For instance, among the polled students, the number of the male students is almost three times than the number of the female students (2.75). Similar unbalanced ratios are also observed between the local and international students (2.55) as well as between the Foundation students and the Undergraduate students (2.12). However, this may imply a similar result for other institutions that possess a larger proportion of male students than female students as well as more local students than international students. Nevertheless, the UNMC boasts of possessing a student body with over 50 different nationalities~(\citealp{unmc12}). This results in a potential for larger sample data of international students as well as further diversity in the sample data. This is a very unique population since none other university in the region possesses such a diverse student body population.

\subsection{Conclusion}

The hypothesized models present the relationships among the level of test anxiety in mathematics, students' expectation score and their actual AP in the context of distinct factors including academic levels, gender groups and nationality backgrounds. After conducting some statistical analyses, it is discovered that the findings are mixed and not all hypotheses are supported. Of the six hypotheses, the findings show that only Hypotheses~\ref{hypothesis1} and \ref{hypothesis3} are supported by analysis, while Hypotheses~\ref{hypothesis2} and Hypothesis~\ref{hypothesis4}--\ref{hypothesis6} are not.

The findings of this research suggests that the level of mathematics-test anxiety during the period of the midterm examination appears to be the most strongly related with the students' expectation.
A negative correlation value with a moderate effect size between the students' expectation score and their level of mathematics-test anxiety is observed ($r = -0.366$, $p = 0.000$). Since the correlation value is significant, this finding supports Hypothesis~\ref{hypothesis1}, but it differs from another finding in the literature, see, for instance~\citeauthor{Chinta05} (\citeyear{Chinta05}). A negative correlation value with a very small effect size is obtained between the students' level of mathematics-test anxiety and their actual AP ($r = -0.099$, $p = 0.157$). This result supports similar findings in the existing literature, for instance those of \citeauthor{Pekrun92}~(\citeyear{Pekrun92}) and \citeauthor{Schwarzer92}~(\citeyear{Schwarzer92}). Since the correlation value is not significant, this result does not support Hypothesis~\ref{hypothesis2}. A positive correlation value with a moderate effect size is discovered between the expectation score and the actual AP ($r = 0.332$, $p = 0.000$). The correlation value is significant and thus it supports Hypothesis~\ref{hypothesis3}.

In the context of academic levels, the students enrolled in the Foundation programme have a higher level of mathematics-test anxiety than those enrolled in the Undergraduate programme. Negative correlation values with very small effect sizes are observed between the students' level of mathematics-test anxiety and their AP for the students enrolled in both the Foundation ($r = -0.103$, $p = 0.200$) and Undergraduate ($r = -0.029$, $p = 0.842$) programmes. Since the correlation values are not significant, this finding does not support Hypothesis~\ref{hypothesis4}. Although there is no particular literature that discusses a similar concern, this result is in line with the findings in other studies on a broader perspective that the older the students are, the weaker the relationship between the level of test anxiety and AP will be, see, for instance \citeauthor{Sud91}~(\citeyear{Sud91}) and \citeauthor{Chapel05}~(\citeyear{Chapel05}).
\begin{table}[h]
{\footnotesize \renewcommand{\tabcolsep}{10pt}
\caption{Summary table} \label{summary}
\begin{tabular}{lclcl}
  \hline \hline
  & Hypothesis & \hspace*{0.5cm} Correlation between & Relevant statistics & \hspace*{0.2cm} Status \\
  \hline
On correlation  & 1 & Students' expectation score and & $r = -0.366$ ($p = 0.000$) & Supported \\
triangle        &   & level of mathematics-test anxiety &  &  \\
                & 2 & Level of mathematics-test   & $r = -0.099$ ($p = 0.157$) & Rejected \\
                &   & anxiety and AP &  &  \\
                & 3 & Students' expectation score  & $r =  0.332$ ($p = 0.000$) & Supported \\
                &   & and actual AP &  &  \\
  &   &  &  &  \\
On AP  & 4 & Foundation and & $r = -0.103$ versus $r = -0.029$ & Rejected \\
       &   &  Undergraduate & ($p = 0.200$) \quad ($p = 0.842$) & \\
       & 5 & Male and female & $r = -0.169$ versus $r = -0.013$ & Rejected  \\
       &   &  & ($p = 0.038$) \quad ($p = 0.925$) & \\
       & 6 & Local Malaysians and   & $r = 0.037$ versus $r = -0.396$ & Rejected \\
       &   & international students & ($p = 0.655$) \quad ($p = 0.002$)   &  \\
 \hline \hline
\end{tabular}
}
\end{table}

In the context of gender difference, female students have a higher level of mathematics-test anxiety than the male students. Further analysis shows that a significant correlation value with a small effect size between the level of mathematics-test anxiety and the AP is only discerned for the male students ($r = -0.169$, $p = 0.038$). For the female students, the correlation value is not significant and has an extremely small effect size ($r = -0.013$, $p = 0.925$). Thus, this finding does not support Hypothesis~\ref{hypothesis5}. Interestingly, this conclusion is contrary to the existing literature since it is discovered that female students have a higher level of test anxiety than male students irrespective of the level of study and cultural background (\citealp{vanderPloeg84,Sharma90,Kleijn94,Karimi09}). Finally, in the context of nationality backgrounds, the local Malaysian students have a higher level of test anxiety in mathematics than their peers from overseas. Further analysis shows that the level of anxiety affects the AP among the international students. A positive correlation value with a very small effect size is observed among the local students, yet it is not significant ($r = 0.037$, $p = 0.655$). On the other hand, a negative, significant correlation value with a moderate effect size is observed among the international students ($r = -0.396$, $p = 0.002$). Thus, although this finding does not support Hypothesis~\ref{hypothesis6}, the result is in agreement with the existing literature (\citealp{Burns91,Stankov10}). A summary table of the hypotheses status (supported and rejected) with relevant statistics is presented in Table~\ref{summary}.

\subsection{Future implication}

The future implications of this study indicate that such a problem of test anxiety does exist among early undergraduate students enrolling in Engineering programme at the UNMC, particularly in subjects which are considered quite challenging for many students, such as Mathematics. Since the set of Likert questionnaire statements is adapted from the TAI (\citealp{Spielberger80}), an identical or a similar set of Likert questionnaire items can be implemented at other tertiary institutions too, whether they are public, government-run institutions or private institutions, including the ones that are run by the local management or branch campuses of overseas universities. Since the institution of this research is very unique, that is, the only branch campus of a British university in Malaysia, it would be interesting to investigate whether similar results can be obtained when this research is conducted in other institutions, whether they are public or private universities. Currently, there are three other universities possessing a similar status of the Malaysian branch campus of overseas universities and all of them are Australian universities. It is most likely that similar results would be obtained if the polled respondents have compositions similar to the ones in this research.

There are a number of techniques that can be used to reduce the level of test anxiety particularly for challenging subjects such as mathematics, physical sciences or engineering modules. Some suggested methods discussed in Subsection~\ref{solution} can also be applied to the students at tertiary institutions, including colleges and universities. Another technique is to control the level of test difficulty, as has been considered by \citeauthor{Hong99} (\citeyear{Hong99}). In addition, the fact that the problem of test anxiety in mathematics exists among the university students gives a clue regarding a possible research direction to improve the mathematics curriculum, in general, and the method of assessment, in particular. There is no doubt that this issue and the related results are worthy of further investigation.

\subsection*{Acknowledgements}
{\small The authors thank Dr. Tom Cross (The University of Nottingham University Park Campus, UK), Linda Ellison (The University of Nottingham University Park Campus, UK), Grace Yap (The University of Nottingham Malaysia Campus), Daniel Hobbs (Sungkyunkwan University, South Korea and Guanghua College, Shanghai Campus, China), Dr. Ardhasena Sopaheluwakan (Indonesian Meteorological, Climatological and Geophysical Agency) and most importantly all anonymous referees for their valuable remarks that have led to the improvement of this article. The authors also thank all the students who participated in this study by completing the questionnaires. \par}

\newpage
\begin{center}
\Large  Appendix\\
\end{center}
\begin{figure}[h]
  \begin{center}
  \includegraphics[width=0.9\textwidth]{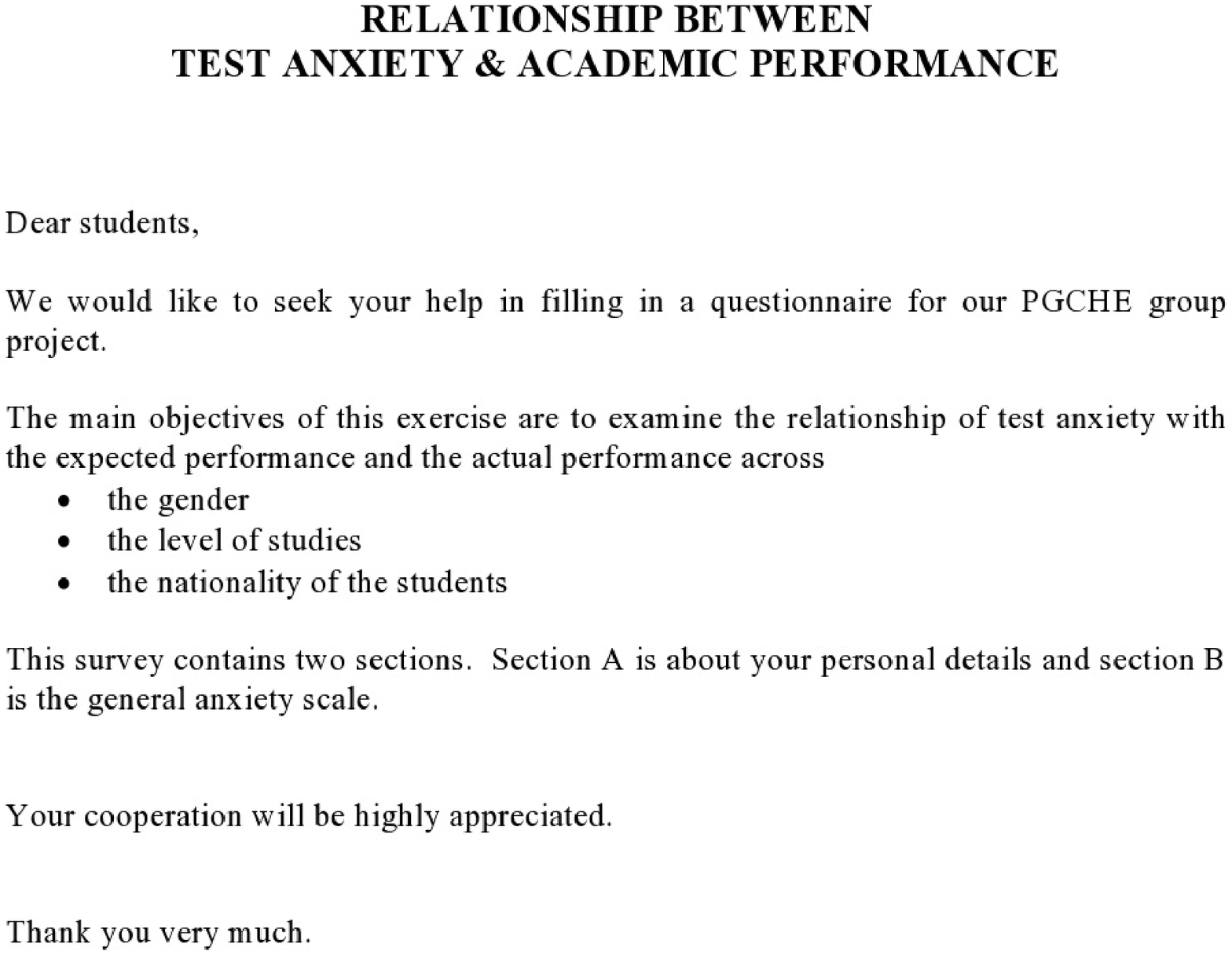}
  \end{center}
\end{figure}

\newpage
\begin{figure}[h!]
  \begin{center}
  \vspace*{1cm}
  \includegraphics[width=0.9\textwidth]{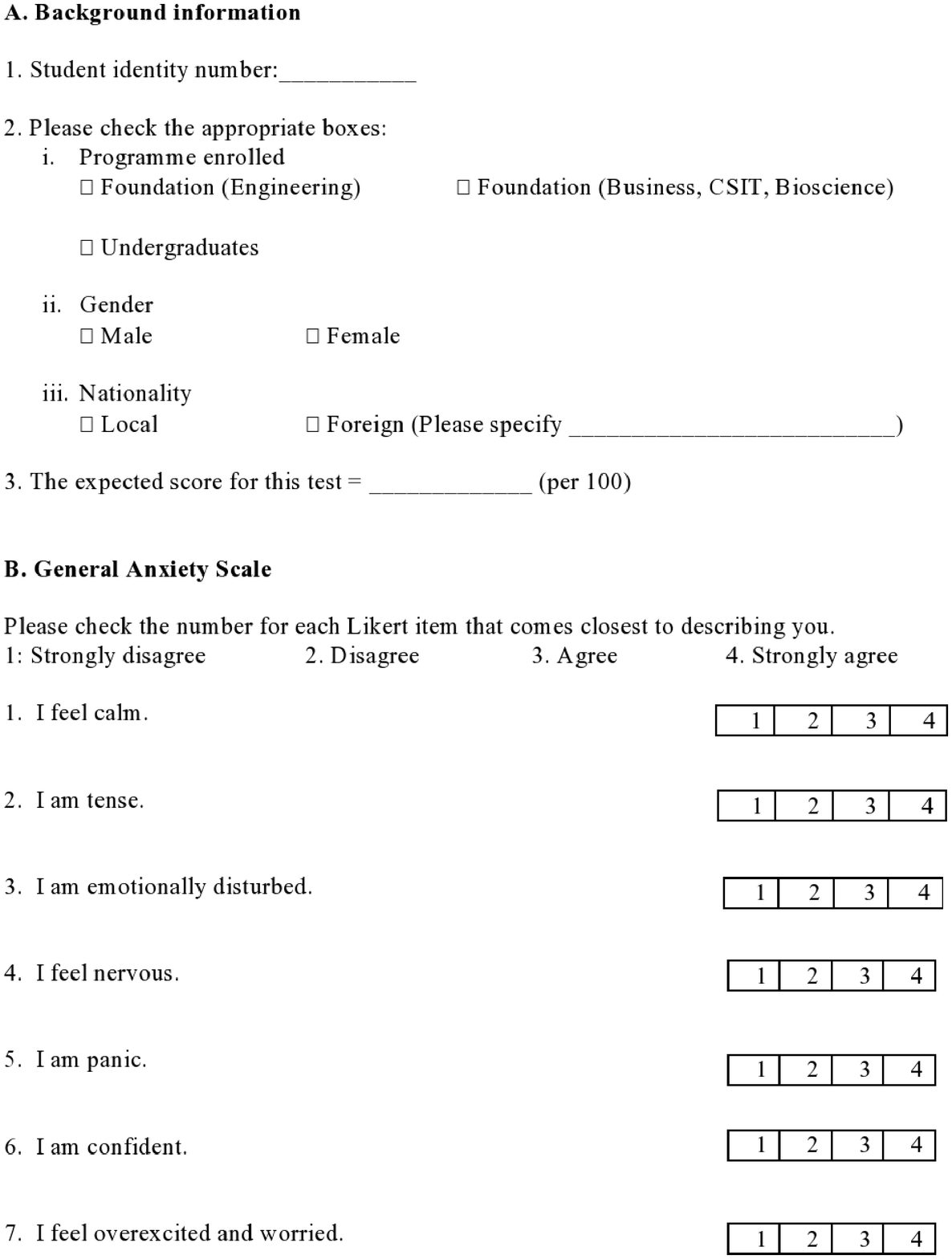}
  \end{center}
\end{figure}
\end{document}